\documentclass[11pt, psamsfonts]{amsart}

\usepackage{a4wide}
\usepackage{amsfonts,amssymb,amscd}
\usepackage{graphicx}
\usepackage{subfigure}
\usepackage{hyperref}

\newtheorem{theorem}{Theorem}[section]

\newtheorem{lemma}[theorem]{Lemma}

\newtheorem{proposition}[theorem]{Proposition}

\theoremstyle{remark}
\newtheorem{conjecture}[theorem]{Conjecture}
\newtheorem{remark}[theorem]{Remark}

\renewcommand{\leq}{\leqslant}
\renewcommand{\geq}{\geqslant}
\renewcommand{\le}{\leqslant}
\renewcommand{\ge}{\geqslant}

\newcommand{\ptl}{\partial}

\newcommand{\Om}{\Omega}

\newcommand{\rr}{\mathbb{R}}

\makeatletter
\renewcommand{\p@subfigure}{}
\makeatother

\numberwithin{equation}{section}

\begin{document}

\title[Least-perimeter partitions of the disk]{Least-perimeter
partitions of the disk into three regions of given areas}

\author[A. Ca\~nete]{Antonio Ca\~nete}
\address{Departamento de Geometr\'{\i}a y Topolog\'{\i}a \\ Facultad
de Ciencias \\ Universidad de Granada \\ E-18071 Granada (Espa\~na)}
\email{antonioc@ugr.es}

\author[M. Ritor\'e]{Manuel Ritor\'e}
\address{Departamento de Geometr\'{\i}a y Topolog\'{\i}a \\ Facultad
de Ciencias \\ Universidad de Granada \\ E-18071 Granada (Espa\~na)}
\email{ritore@ugr.es}

\thanks{Both authors have been supported by MCyT-Feder research
project BFM2001-3489}
\date{June 25, 2003}


\begin{abstract}
We prove that the unique least-perimeter way of partitioning the unit
$2$-dimen\-sio\-nal disk into three regions of prescribed areas is by
means of the standard graph described in Figure~\ref{fig:standard}.
\end{abstract}

\maketitle

\thispagestyle{empty}

\section*{Introduction}
\label{sec:introduccion}

Partitioning problems in the Calculus of Variations have multiple
applications in physical sciences.  They can model multitude of
natural phenomena such as the shape of a cellular tissue, the
interface of separation between fluids, and many others, as described
in the treatise by D'Arcy Thompson~\cite{thompson}.

\begin{figure}[h]
\centerline{\includegraphics[width=0.2\textwidth]{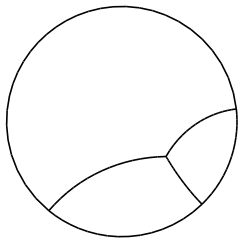}}
\caption{The least-perimeter partition of the disk into three given
areas}
\label{fig:standard}
\end{figure}

In this work we consider the isoperimetric problem of partitioning a
planar disk into three regions of given areas with the least possible
perimeter, and we prove that the standard configuration in
Figure~\ref{fig:standard}, consisting of three circular arcs or
segments meeting orthogonally the boundary of the disk, and meeting in
threes at $120$ degrees in an interior vertex, is the only solution to
this problem.

In addition to the above conditions, the solution must satisfy a
certain balancing condition on the geodesic curvatures of the circles.
This condition will be stated precisely in the next section.

Existence and regularity of solutions for this problem are guaranteed
by the results of F.~Morgan~\cite{morgan-soap}, who showed that the
minimizer, in the interior of the disk, is composed of smooth curves
of constant geodesic curvature meeting in threes at $120$ degree
angles.  Boundary regularity also follows from \cite{morgan-soap}
although it is not explicitly stated in his work.  Existence and
regularity in higher dimension were studied by
F.~Almgren~\cite{almgren}. 

The least-perimeter way of partitioning a disk $D$ into two regions of
given areas is by means of an arc of circle or segment that meets
orthogonally $\ptl D$.  From the existence and regularity results in
next section it follows that there is a solution, which is a smooth,
possibly nonconnected, embedded curve with constant geodesic curvature
that meets $\ptl D$ orthogonally.  Such a curve must be connected,
since otherwise we could rotate one component with respect to the
center of the disk until it touches a second one, thus producing a
non-allowed singularity.  On the other hand, as the curve has constant
geodesic curvature, it must be part of a circle or of a line.

\begin{figure}[h]
\centerline{\includegraphics[width=0.2\textwidth]{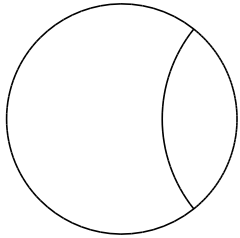}}
\caption{The least-perimeter partition of the disk into two given
areas}
\label{fig:2-standard}
\end{figure}

The isoperimetric problem consisting of enclosing $n$ given areas in
the disk or in the plane with the least possible perimeter has a
strong complexity which is derived not from the geometry of the
individual components of the solution (they can be described in terms of
circles or lines) but from their large number.

The planar double bubble conjecture was proved by J.~Foisy et al. 
\cite{foisy-zimba}, who showed in 1993 that the standard planar double
bubble uniquely minimizes perimeter in $\rr^2$.  Assuming that the
studied regions are connected, C.~Cox et al.~\cite{C} proved in 1994
that the standard planar triple bubble uniquely minimizes perimeter in
the plane for any three given areas.  R.~P.~Devereaux~\cite{V} studied
in 1998 the planar triple bubble conjecture under the hypothesis that
all the regions have the same pressure.  W.~Wichiramala finally proved
the planar triple bubble conjecture in 2002 in his
Ph.~D.~Thesis~\cite{W}.  J.~Masters~\cite{M} proved in 1996 the double
bubble conjecture in $\mathbb{S}^2$.  Interesting preliminary work was
carried out by Bleicher \cite{bleicher-1}, \cite{bleicher-2},
\cite{bleicher-3}.  Concerning boundary problems, G.~Hruska et
al.~\cite{hruska} have obtained some results for planar bubbles in
corners.  Also results on tori and cones have been obtained in
\cite{corneli}, \cite{borawski}.

In higher dimensions, J.~Hass and R.~Schlafly \cite{hs} proved the
double bubble conjecture in $\rr^3$ for equal volumes.  The general
conjecture was settled by M.~Hutchings et al.~\cite{HMRR}.  For higher
dimensional Euclidean spaces Reichardt et al.~\cite{reich} have
obtained a proof of the double bubble conjecture in $\rr^4$ and
partial results in higher dimensional Euclidean spaces.  In the
three-dimensional torus, M.~Carrion et al.~\cite{carrion} have
provided numerical evidence for a double bubble conjecture with ten
types of solutions.  In the three-dimensional sphere and the
three-dimensional hyperbolic space, A.~Cotton and D.~Freeman~\cite{CF}
have also obtained partial results on the conjecture that the standard
double bubble in these spaces uniquely minimizes perimeter.

Planar bubbles are also of great interest to physicists.  Interesting
articles focusing on physical aspects of the problem are
\cite{graner-1} and \cite{graner-2}.

The most interesting mathematical open question for these problems is
to show that the minimizing configurations must have connected
regions, either in $\rr^2$ or in the disk.  In addition, in the planar
problem, one should also be able to prove that the exterior region is
connected, i.e., that there are no empty chambers.

We have organized this paper in several sections.  In
Section~\ref{sec:preliminare} we give precise definitions, compute the
first and second variations of length for graphs, recall existence and
regularity results for the problem of minimizing perimeter while
partitioning the disk into given areas, and state some properties
which minimizing graphs must satisfy.  In Section~\ref{sec:cota} we
obtain a bound on the number of components of the largest pressure
region determined by a graph which minimizes perimeter up to second
order.  We conclude that a minimizing configuration must have one of
ten possible types, described in Figure~\ref{fig:configs}.  In
Section~\ref{sec:unstable} we prove the necessary results to discard
the~possibilities obtained in Section~\ref{sec:cota}, which allow us
to prove our Main Theorem in Section~\ref{sec:teorema}.  In a final
section, we indicate further lines of research and give several
conjectures.

All the pictures in this paper have been made by using Surface
Evolver, a software developed by Ken Brakke
(http://www.susqu.edu/facstaff/b/brakke/).

\section{Preliminaries}
\label{sec:preliminare}

\subsection{Notation}
Let $D\subset\rr^2$ be the closed unit disk in centered at the origin. 
An {\em admissible graph} $C\subset D$ consists of vertices and curves
so that at every interior vertex (that is, a vertex in the interior of $D$) three
curves of $C$ meet and at every boundary vertex (a vertex in $\ptl D$) just
one curve of $C$ meets $\ptl D$.  We shall also assume that $C$
induces a decomposition of the open unit disk into $n$ regions
$R_{i}$, $1\le i\le n$, possibly nonconnected.  An {\em $m$-component}
is a connected component of a region with $m$ edges.

If $R_{i}$ and $R_{j}$ are adjacent regions, we will denote by
$C_{ij}\subset C$ the (not necessarily connected) curve separating them. 
Let $I(i)=\{j\neq i; R_{j}\ \text{touches}\ R_{i}\}$.  With this
notation
\[
\ptl R_{i}\cap\text{int}(D)=\bigcup_{j\in I(i)} C_{ij}.
\]
We shall denote by $N_{ij}$ the normal vector to the curve $C_{ij}$ 
pointing into the region $R_{i}$, and by $h_{ij}$ the geodesic 
curvature of the curve $C_{ij}$ with respect to the normal $N_{ij}$.

A {\em standard graph} consists in three circular arcs or lines
segments meeting at an interior vertex at 120 degree angles, reaching
orthogonally $\partial D$, and so that the sum of the geodesic
curvatures is zero.

Given $n$ positive numbers $a_1,\ldots,a_n$ such that $\sum_{i=1}^n
a_{i}=\pi$, the {\em isoperimetric profile} is the function
$I(a_{1},\ldots,a_{n})$ defined as the infimum of the lengths of all
admissible graphs separating regions in the disk of areas
$a_1,\ldots,a_n$.

We will say that an admissible graph $C$ is {\em minimizing} for
prescribed areas $a_1,\ldots,a_n$ if $I(a_1, \ldots, a_n)$ is attained
by $C$.

\subsection{Variational formulae}
Given an admissible graph $C\subset D$, we will consider smooth
one-parameter variations $\varphi_t:C\to D$ for $t$ small, which
satisfy $\varphi_t(\ptl D)\subset\ptl D$.  We will denote by $X=d
\varphi_t/dt|_{t=0}$ the associated infinitesimal vector field, which
is smooth on every curve $C_{ij}$.  Note that $X(p)$ is tangent to
$\ptl D$ for each $p$ in $\ptl D$.  Let $u_{ij}=X\cdot N_{ij}$ be the
normal component of $X$ on $C_{ij}$.

Given such a variation, it is easy to check that the derivative of the
area $A_{i}$ of $R_i$ at $t=0$ is given by
\begin{equation}
\label{eq:dareai}
\frac{dA_{i}}{dt}\bigg|_{t=0}=-\sum_{j\in I(i)} \int_{C_{ij}} u_{ij}.
\end{equation}

For the derivative of length for such a variation we have

\begin{proposition}[First variation of length {\cite[Lemma~3.1]{HMRR}}]
\label{prop:firstvariation}
Consider an admissible graph $C\subset D$, and a smooth variation
$\varphi_{t}:C\to D$ with associated vector field $X$.  Then the first
derivative of the length of $\varphi_{t}(C)$ at $t=0$ is given by
\begin{equation}
\label{eq:firstvariation}
\frac{dL}{dt}\bigg|_{t=0}=-\frac{1}{2}\,\sum_{\substack{i\in\{1,\ldots,n\}
\\j\in I(i)}} \bigg\{\int_{C_{ij}} h_{ij}u_{ij}+\sum_{p\in\partial
C_{ij}} X(p)\cdot\nu_{ij}(p)\bigg\},
\end{equation}
where $\nu_{ij}(p)$ is the inner conormal to $C_{ij}$ in $p$ .
\end{proposition}

We will say that an admissible graph is {\em stationary} if
\eqref{eq:firstvariation} vanishes for any area-preserving variation. 
From Proposition~\ref{prop:firstvariation} it is easy to prove the following

\begin{proposition}
\label{prop:conditions}
Let $C\subset D$ be a stationary graph.  Then the following conditions
are satisfied
\begin{itemize}
\item[(i)] The geodesic curvature $h_{ij}$ is constant on $C_{ij}$.

\item[(ii)] The edges of $C$ meet in threes at $120$-degree
angles in interior vertices.  

\item[(iii)] The balancing condition: three edges $C_{ij}$, $C_{jk}$,
$C_{ki}$ meeting in an interior vertex satisfy\begin{equation}
\label{eq:cocycle}
h_{ij}+h_{jk}+h_{ki}=0.
\end{equation}

\item[(iv)] The edges of $C$ meet $\partial D$ orthogonally at 
boundary vertices.
\end{itemize} 
\end{proposition}

Condition (ii) implies that, in some interior vertex where the three
curves $C_{ij}$, $C_{jk}$, $C_{ki}$ meet, the normals add up to zero,
i.~e., $N_{ij}+N_{jk}+N_{ki}=0$.  This implies that the normal
components of the vector field $X$ must satisfy
\begin{equation}
\label{eq:uijk}
u_{ij}+u_{jk}+u_{ki}=0.
\end{equation}

Given a stationary graph $C$, and a function $u:\bigcup_{i,j}
C_{ij}\to\rr$, with $u_{ij}=u|_{C_{ij}}$, satisfying condition
\eqref{eq:uijk} on every interior vertex, it is always possible to
find a vector field $X$ on $C$, so that $u_{ij}=X\cdot N_{ij}$ and $X$
is tangent to $\ptl D$ in each boundary vertex.  Associated to $X$ one
can also find a one-parameter variation $\varphi_{t}:C\to D$, for $t$
small enough, so that $\varphi_{t}(p)=\exp_{p}(tX(p))$ for any $p$ out
of an arbitrarily small neighbourhood of $\ptl D$.  The argument is as
follows: fix an arbitrary neighbourhood $U$ of $\ptl D$ that does not
contain interior vertices of $C$.  Modify $X$ so that it is normal to
$C$ in $U$.  Let $\nu$ be the inner normal to $\ptl D$.  Extend it to
$U$ so that it is tangent to the edges of $C$.  Also extend $X$ to a
vector field on $U\cap D$ by means of the exponential mapping.  Let
$\lambda$ be a smooth function equal to $1$ near $\ptl D$ with support
in $U$.  Consider the vector field
$Y=X-(X\cdot(\lambda\nu))\,\lambda\nu$ and the local one-parameter
group $\psi_{t}$ generated by $Y$.  Since $\nu$ is tangent to $C$ and
$X$ is normal to $C$ in $U$, we have that $Y=X$ on $C$.  Moreover, for
$p\in\ptl D$, the vector $Y(p)$ is tangent to $\ptl D$.  Hence the
deformation $\psi_{t}(C\cap U)$ has initial velocity vector field $X$
and keeps $C$ inside the disk.  The variation $\psi_{t}(C)$ has the
further property that coincides with $\exp_{p}(tX(p))$ in $U$ out of
the support of $\lambda$.  Now we simply define $\varphi_{t}(p)$ equal
to $\exp_{p}(tX(p))$ out of the support of $\lambda$, and equal to
$\psi_{t}(p)$ in $U$.

The balancing condition \eqref{eq:cocycle} allows us to define 
a pressure $p_{i}$ on every region $R_{i}$, starting from a given 
region, so that 
\begin{equation}
\label{eq:pressures}
h_{ij}=p_{i}-p_{j}.
\end{equation}
These pressures are determined up to an additive constant. The first 
variation formula of length can be rewritten in terms of pressures 
in the following way: if $C$ is a stationary graph, then the first 
variation of length for an arbitrary variation is given by:
\begin{equation}
\label{eq:1stpressures}
\frac{dL}{dt}=\sum_{i=1}^n p_{i}\,\frac{dA_{i}}{dt}.
\end{equation}
Observe that the indetermination of the pressures up to some additive 
constant does not affect the above formula since $\sum_{i=1}^n 
dA_{i}/dt=0$ for any variation of the regions $R_{i}$, as $\sum_{i=1}^n
A_{i}(t)=\text{area}(D)$ along the variation.

Let us prove now the second variation formula of length

\begin{proposition}[Second variation of length]
\label{prop:secondvariation}
Let $C$ be a stationary graph and let $\{\varphi_{t}\}$ be a variation
with associated vector field $X$ preserving areas up to second order. 
Then the second derivative of length at $t=0$ is given by
\begin{align}
\label{eq:secondvariation}
-\frac{1}{2}\,\sum_{\substack{i=1,\ldots,n \\ j\in I(i)}} \bigg\{
\int_{C_{ij}}(u_{ij}'' +h_{ij}^2 u_{ij})\,u_{ij}
+\sum_{\substack{p\in\partial C_{ij} \\p\in\text{\em int}(D)}}
\bigg(-q_{ij}u_{ij}^2 &+ u_{ij}\frac{\partial u_{ij}}{\partial
\nu_{ij}} \bigg)(p) \\
\nonumber &+ \sum_{\substack{p\in\partial C_{ij} \\p\in\partial
D}}\bigg(u_{ij}^2+u_{ij}\frac{\partial u_{ij}}{\partial
\nu_{ij}}\bigg)(p) \bigg\},
\end{align}
where $q_{ij}(p)=(h_{ki}+h_{kj})(p)/\sqrt{3}$, and $R_{k}$ is the third 
region touching the vertex $p$.
\end{proposition}
\begin{proof}
Differentiating the integral terms in
equation~\eqref{eq:firstvariation}, we get
\[
\frac{d}{dt}\bigg|_{t=0} \bigg(\int_{C_{ij}}\,h_{ij}u_{ij}\bigg)=
\int_{C_{ij}}(u_{ij}'' +h_{ij}^2 u_{ij})\,u_{ij}
+h_{ij}\,\frac{d}{dt}\bigg |_{t=0}\bigg(\int_{C_{ij}} u_{ij}\bigg),
\]
but since the variation preserves areas up to second order, it follows
that
\[
\sum_{\substack{i\in\{1,\ldots,n\} \\ j\in I(i)}}
h_{ij}\,\frac{d}{dt}\bigg|_{t=0}\bigg(\int_{C_{ij}} 
u_{ij}\bigg)=2\,\sum_{i=1}^n p_{i}\,\frac{d^2\!A_{i}}{dt^2}\bigg|_{t=0}=0.
\]
Differentiating now the second term in
equation~\eqref{eq:firstvariation}, we get
\[
\frac{d}{dt}\bigg|_{t=0} (X\cdot\nu_{ij}) = (D_X X\cdot\nu_{ij}) +
u_{ij}\,h_{ij}\, (X\cdot\nu_{ij}) + u_{ij}\frac{\partial
u_{ij}}{\partial\nu_{ij}}.
\]
For $p\in\text{int}(D)$, the first term vanishes since
$\nu_{ij}+\nu_{jk}+\nu_{ki}=0$, and after some calculations as in
\cite{HMRR}, the second one can be seen as $-q_{ij}u_{ij}^2$, where
$q_{ij}=(h_{ki}+h_{kj})/\sqrt{3}$.  For $p\in\partial D$, since the
configuration is stationary, the edges meet $\partial D$ orthogonally,
so that $D_X X(p)\cdot\nu_{ij}(p)$ equals $u_{ij}^2$ times the
geodesic curvature of $\partial D$, and $(X\cdot\nu_{ij})(p)=0$.
\end{proof}

The condition that the variation must preserve area up to second order 
is not really needed as we can show in the next Lemma

\begin{lemma}
\label{lem:ap}
Let $C\subset D$ be a stationary graph.  Given smooth functions
$u_{ij}:C_{ij}\to\rr$ such that \eqref{eq:dareai} and \eqref{eq:uijk}
are satisfied $($a variation that preserves area up to first order is
given$)$, there is a variation $\{\varphi_{t}\}$ of $C$ which leaves
constant the area of the regions enclosed by $\varphi_{t}(C)$ and such
that the normal components of the initial velocity vector field $X$
are the functions~$u_{ij}$.
\end{lemma}

\begin{proof}
Let $X$ be a vector field on $C$, smooth over each curve $C_{ij}$,
such that $X\cdot N_{ij}=u_{ij}$.  Let $\psi_t:C\to D$ be a
one-parameter variation of $C$ associated to $X$ such that
$\psi_{t}(p)=\exp_{p}(tX(p))$ out of a small neighbourhood $U$ of
$\ptl D$ which does not contain interior vertices of $C$.

We label the regions $R_{i}$ so that $R_{i}$ touches $R_{i+1}$ for 
$i=1$, $\ldots$, $n-1$. Choose positive functions $v_{i}$ with
support in the interior of $C_{i(i+1)}$ and out of $U$. The variation 
induced by the vector field $v_{i}\,N_{i(i+1)}$ decreases the area of 
$R_{i}$, increases the area of $R_{i+1}$ and leaves constant the area 
of the remaining regions.

Consider the variation equal to
\[
(t,s_{1},\dots,s_{n-1})\longmapsto \exp_{p}\big(tX(p)
+\sum_{i=1}^{n-1} s_{i}v_{i}N_{i(i+1)}(p)\big), \qquad\text{in }C\cap(D-U),
\]
and equal to $\psi_{t}(p)$ for $p\in C\cap U$. Consider the function 
$(A_{1},\dots,A_{n-1})$ of $(t,s_{1},\ldots,s_{n-1})$, given by the 
areas of the deformation of the regions $R_{1}$,$\ldots$, $R_{n-1}$. 
The Jacobian
\[
\frac{\ptl(A_{1},\ldots,A_{n-1})}{\ptl(s_{1},\ldots,s_{n-1})}
\]
is lower triangular, with non-vanishing entries in the principal diagonal, 
so that the matrix is regular. The Implicit Function Theorem allows 
us to find smooth functions $s_{1}(t)$,$\ldots$, $s_{n-1}(t)$ such 
that $A_{i}(t,s_{1}(t),\ldots,s_{n-1}(t))$ is constant for all $i$.

The initial velocity vector field of such a variation is equal to $X$
on $C\cap U$, and to $X+\sum_{i=1}^{n-1} s_{i}'(0)\,v_{i}N_{i(i+1)}$
on $C\cap(D-U)$.  As $s_{i}'(0)=0$ since $\psi_{t}$ preserves areas up
to first order, we conclude that $X$ is the initial velocity vector
field.
\end{proof}

\begin{remark}
\label{rem:stationary}
A variation of a stationary graph $C$ by stationary graphs preserves
the angles between edges at interior vertices and the orthogonality
condition at boundary vertices.  Given a variation preserving the area
of all the regions up to first order, we can modify it by
Lemma~\ref{lem:ap} so that the areas enclosed are constant along the
deformation.  From the second variation formula we get that the second
derivative of length is given by
\[
\frac{d^2\!L}{dt^2}=\sum_{\alpha} 
\frac{dp_{\alpha}}{dt}\,\frac{dA_{\alpha}}{dt},
\]
where $\alpha$ labels the {\em components} of the stationary graph
(regions can be disconnected), and $dp_{\alpha}/dt$ is the derivative
of the pressure of the component $\alpha$ with respect to the
considered variation.  Take into account that the quantity
$u_{ij}''+h_{ij}^2 u_{ij}$, the derivative of the geodesic curvature
$h_{ij}$, only depends on $u_{ij}$, the normal component of the
variational vector field $X$, and that the modification needed in
Lemma~\ref{lem:ap} to preserve areas only modifies the acceleration of
the variation.  The angle-preserving condition depends only on the
initial velocity vector field.
\end{remark}

In general, if the areas are not preserved up to second order, the
second derivative of length, for a deformation of a stationary graph
by stationary graphs, is given by
\[
\frac{d^2\!L}{dt^2}=\sum_{\alpha}
\frac{dp_{\alpha}}{dt}\,\frac{dA_{\alpha}}{dt}
+p_{\alpha}\frac{d^2\!A_{\alpha}}{dt^2},
\]
which can also be obtained by differentiating
equation~\eqref{eq:1stpressures}.

\begin{remark}
For a variation such that the angles between the edges are preserved,
we have $D_X (\nu_{ij}+\nu_{jk}+\nu_{ki})=0$ (since
$\nu_{ij}+\nu_{jk}+\nu_{ki}=0$ for all $t$), so the boundary term in
the second variation formula vanishes.
\end{remark}

\subsection{Admissible functions and the index form}
Let $C$ be a stationary graph.  We say a function $u:\bigcup_{i,j}
C_{ij}\to\rr$ is {\em admissible} if the restrictions
$u_{ij}=u|_{C_{ij}}$ lie in the Sobolev space $W^{1,2}(C_{ij})$, and
verify that at any interior vertex $p$, $u_{ij}(p) + u_{jk}(p) +
u_{ki}(p) =0$.  These functions correspond to variations of $C$ which
have as normal components of the associated vector field $X$ the
functions $u_{ij}$.  These variations will preserve areas if, for each
$i$,
\[
\sum_{j\in I(i)} \int_{C_{ij}} u_{ij} = 0.
\]

An admissible function $u$ is said to be a {\em Jacobi function} if
the associated variation preserves the geodesic curvatures of each
edge $C_{ij}$ and the angles in each vertex.  The fact that the
geodesic curvatures are preserved means that the restrictions $u_{ij}$
to $C_{ij}$ verify
\[
u_{ij}'' + h_{ij}^2 u_{ij}=0.
\]

It is clear that the normal component of the Killing vector field
generated by the rotations about the origin gives a Jacobi function.

From equation~\eqref{eq:secondvariation}, we define the {\em index
form}, that is, the induced bilinear form defined on the space of
admissible functions, by
\begin{align}
\label{eq:indexform}
Q(u,v)=-\frac{1}{2}\,\bigg\{&\sum_{\substack{i=1,\ldots,n \\ j\in I(i)}} 
\int_{C_{ij}}(u_{ij}'' +h_{ij}^2 u_{ij})\,v_{ij}
\\
\nonumber
+&\sum_{\substack{p\in\partial C_{ij} \\p\in\text{int}(D)}}
\bigg(-q_{ij}u_{ij} + \frac{\partial u_{ij}}{\partial
\nu_{ij}} \bigg)(p)\,v_{ij}(p)
+ \sum_{\substack{p\in\partial C_{ij} \\p\in\partial
D}}\bigg(u_{ij} + \frac{\partial u_{ij}}{\partial
\nu_{ij}}\bigg)(p)\,v_{ij}(p) \bigg\},
\end{align}
where $q_{ij}$ are the functions defined in
Proposition~\ref{prop:secondvariation}.

We will say a stationary graph $C$ is {\em stable} if
$Q(u,u)\geq 0$ for any admissible function $u$ whose associated
variation preserves areas, and {\em unstable} if it is not stable.  It
is clear that a minimizing configuration must be stable.

\subsection{Existence and Regularity}

From the results of F.~Morgan~\cite{morgan-soap}, we obtain the
fo\-llowing

\begin{theorem}[{Existence and Regularity~\cite[Th.~2.3]{morgan-soap}}]
\label{te:exist}
Let $D\subset \rr^2$ be a closed disk, and let $a_1,\ldots,a_n$ be $n$
given areas such that $\sum_{i=1}^n a_i=\text{\em area}(D)$.  Then
there exists a graph separating $D$ into $n$ regions of areas
$a_1,\ldots,a_n$.  Moreover such a graph consists of constant geodesic
curvature curves meeting in threes in the interior of $D$ at $120$
degree angles, satisfying the balancing condition \eqref{eq:cocycle}
for the geodesic curvatures, and meeting $\ptl D$, one at a time, in
an orthogonal way.
\end{theorem} 

\begin{proof}
From the results in \cite{morgan-soap} one gets the existence of a
solution and the regularity in the interior of the disk with just
triple points as possible singularities.  One also gets that there is
a finite number of components (and hence of curves) in the minimizing
configuration.  For the boundary regularity, we only need to prove
that at every point of $\ptl D$, at most one curve of the minimizing
configuration arrives, at $90$ degrees.

If one or several curves meet $\ptl D$ at $p$ and at least one of them
is not orthogonal to $\ptl D$, then the first variation formula
implies that the graph is not stationary.  Suppose now that several
curves meet orthogonally $\ptl D$ at $p$.  We order them
counter-clockwise and we consider the first one, $C$, which is the
common boundary of components $\Om_i$ and $\Om_j$, with $\Om_i$ a
boundary one.  Make a small deformation in the interior of $C$ which
implies a loss of area $\delta$ for $\Om_i$.  In order to preserve the
areas, it is possible to choose a point $q$ near $p$ in $C$, join $q$
to the second curve $C'$, which is in the boundary of $\Om_{j}$, and
eliminate the part of $C$ between $p$ and $q$.  It can be checked that
this new configuration, for $\delta$ small enough, reduces perimeter. 
Then we get the desired regularity in the boundary of $D$.
\end{proof}

\subsection{Some properties of minimizing graphs}
\label{sec:geopro}
We now give and recall some results on minimizing graphs that will 
be used to prove our main theorem

\begin{lemma}
\label{le:cotaperfil}
Given $n$ positive numbers $a_{1} ,\ldots, a_{n}$ such that 
$\sum_{i=1}^n a_{i}=\pi$, we have 
\begin{equation}
\label{eq:ene}
I(a_{1}, \ldots, a_{n})\leq n.
\end{equation}
Moreover, equality is never achieved for $n\ge 4$.  If equality holds
in the case $n=3$ then the standard graph consisting of three line
segments dividing the disk into three regions of equal areas is
minimizing.
\end{lemma}

\begin{proof}
We can divide the disk into regions of given areas $a_{1}, \ldots,
a_{n}$ by using appropriate $n$ radii.  This gives \eqref{eq:ene}. 
For $n\ge 4$, this configuration has a prohibited singularity at the
origin, so that it cannot be minimizing.  If equality holds in
\eqref{eq:ene} for $n=3$, the configuration must be stationary, so
that the three radii meet in 120 degrees, and the configuration is the
standard one for equal areas.
\end{proof}

\begin{lemma}
\label{le:twice}
A minimizing graph must be connected.
\end{lemma}

\begin{proof}
On a nonconnected graph, we can rotate one of the components until it
touches another one creating an irregular meeting, so the graph
cannot be minimizing.
\end{proof}

\begin{remark}
Let $C\subset D$ be a minimizing graph, and $\Om$ a connected 
component of $D-C$. Lemma~\ref{le:twice} implies that $\ptl \Om \cap 
\ptl D$ has to be connected.
\end{remark}

\begin{lemma}[{\cite[Lemma~2.4]{foisy-zimba}}]
\label{le:2-component}
On a minimizing graph, there are no $2$-components.
\end{lemma}

\section{A bound on the number of components of the largest pressure 
region}
\label{sec:cota}

\begin{lemma}
\label{le:cota}
Let $C$ be a stable graph separating $D$ into $n$ regions.
Then the region of largest pressure has at most $n-1$ nonhexagonal 
components.
\end{lemma}

\begin{proof}
Assume $R_1$ is the region of largest pressure and suppose it has at
least $n$ nonhexagonal components, $\Om_1, \ldots, \Om_n$.  For each $i$,
consider the variation given by $u_i = 1$ on $\partial \Om_i$, 
extended by zero to the whole graph. If $\Om_{i}$ is a boundary 
component then 
\[
Q(u_{i},u_{i})= 
-\sum_{j\in I(1)}\bigg\{
\int_{C_{1j}\cap\ptl\Om_{i}}h_{1j}^2 
+\sum_{\substack{p\in\partial C_{1j}\cap\ptl\Om_{i} 
\\p\in\text{int}(D)}}
-q_{1j}(p)
+ \sum_{\substack{p\in\partial C_{1j}\cap\ptl\Om_{i} \\p\in\partial
D}} 1\bigg\}
<0,
\]
since, for $p$ in $C_{1j}\cap C_{1k}\cap\ptl \Om_i$, we have
\[
q_{1j}(p)+q_{1k}(p)=\frac{h_{k1}+h_{kj}+h_{j1}+h_{jk}}{\sqrt{3}}(p)
=\frac{h_{k1}+h_{j1}}{\sqrt{3}}(p)\leq 0,
\]
as $R_1$ has the largest pressure.

If $\Om_{i}$ is an interior component then $Q(u_{i},u_{i})$ can be
computed as above except that the last summand does not appear.  So we
get
\[
Q(u_{i},u_{i})\leq 0,
\]
and equality holds if and only if $\Om_i$ is bounded by segments.  It
is easy to obtain, from Gauss-Bonnet Theorem, that $\Om_i$ has to be
an hexagon.  In the case of three regions, this only occurs if the
three pressures are equal.

Hence, in our case we can find some nontrivial
linear combination $u$ of $u_i$, such that the induced variation
preserves areas up to first order and $Q(u,u)<0$.  
\end{proof}

\begin{lemma}
\label{le:possibilities}
Let $C\subset D$ be a minimizing graph separating $D$ into three
regions.  Then $C$ is one of the graphs in Figure~\ref{fig:configs}.
\end{lemma}

\begin{proof}
Suppose first that all the pressures are equal.  If all the components
touch the boundary of $D$ then $C$ is standard.  If there is an
interior component, then it is hexagonal.  It is easy to see that the
edges leaving the vertices of the hexagon meet $\ptl D$ (otherwise we
could find two different parallel rays meeting orthogonally $\ptl D$). 
This implies that the graph is like in Figure~\ref{fig:hex}.  This
graph has two regions with three nonhexagonal convex components, and
so it is unstable by Lemma~\ref{le:cota}.

\begin{figure}[h]
\centering{\includegraphics[width=0.2\textwidth]{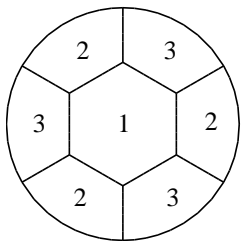}}
\caption{A graph with an hexagonal interior component}
\label{fig:hex}
\end{figure}

\begin{figure}[htp]
\centering{
\subfigure[$\alpha$,
$\beta\in\{1,2\}$]{\label{conf1}\includegraphics[width=0.2\textwidth]{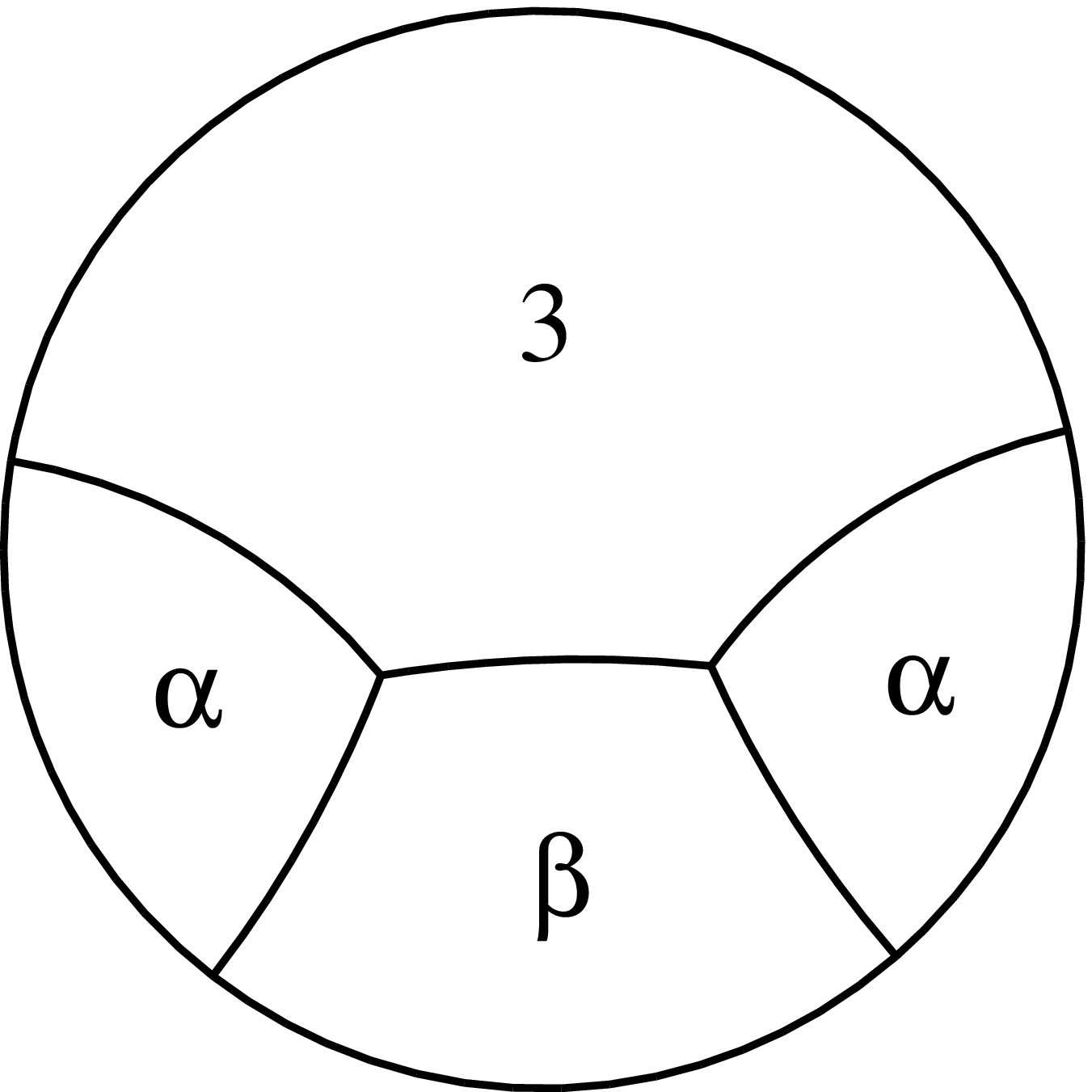}}
\hspace{0.1\textwidth}
\subfigure[]{\label{conf3}\includegraphics[width=0.2\textwidth]{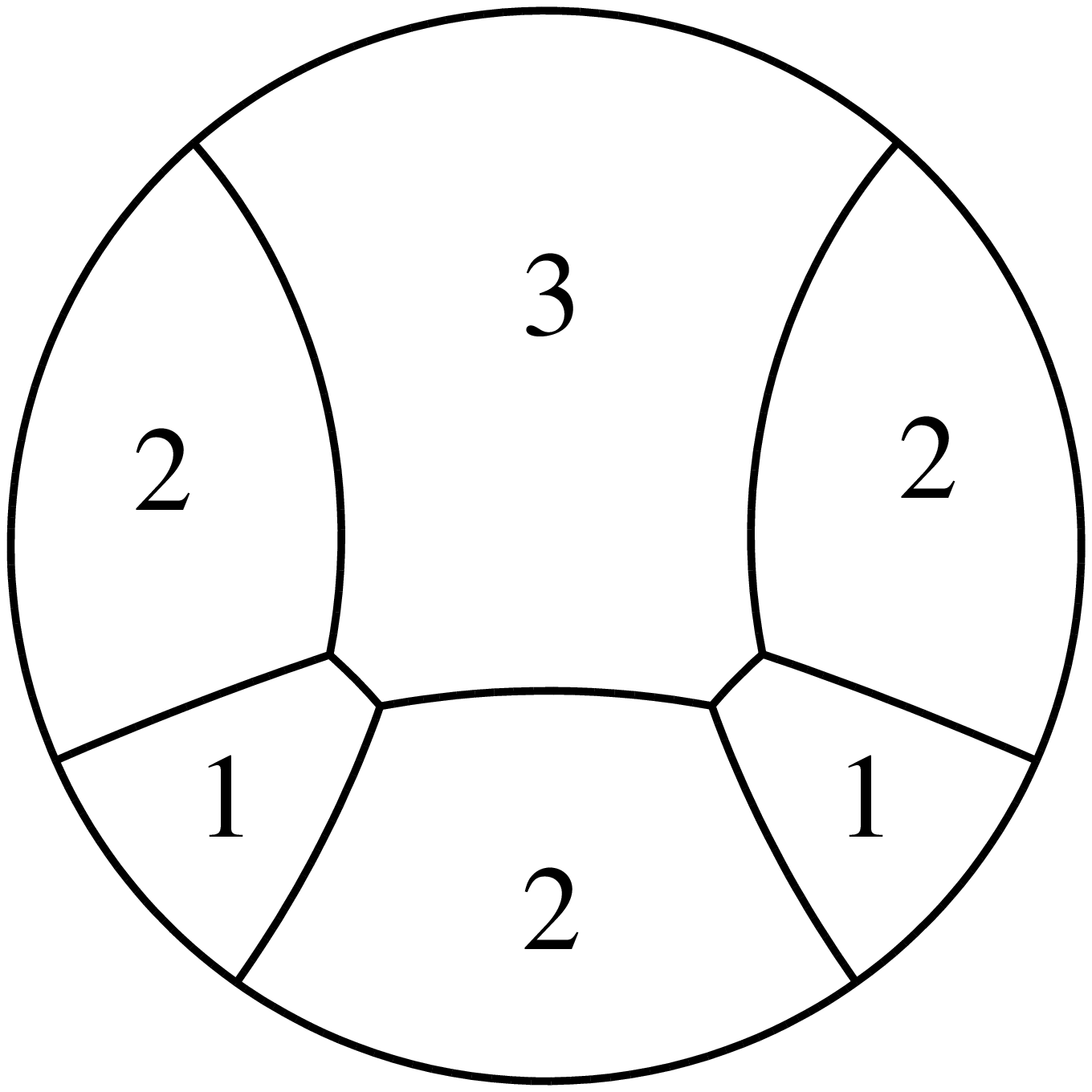}}
\hspace{0.1\textwidth}
\subfigure[$\alpha$,
$\beta\in\{1,2\}$]{\label{conf4}\includegraphics[width=0.2\textwidth]{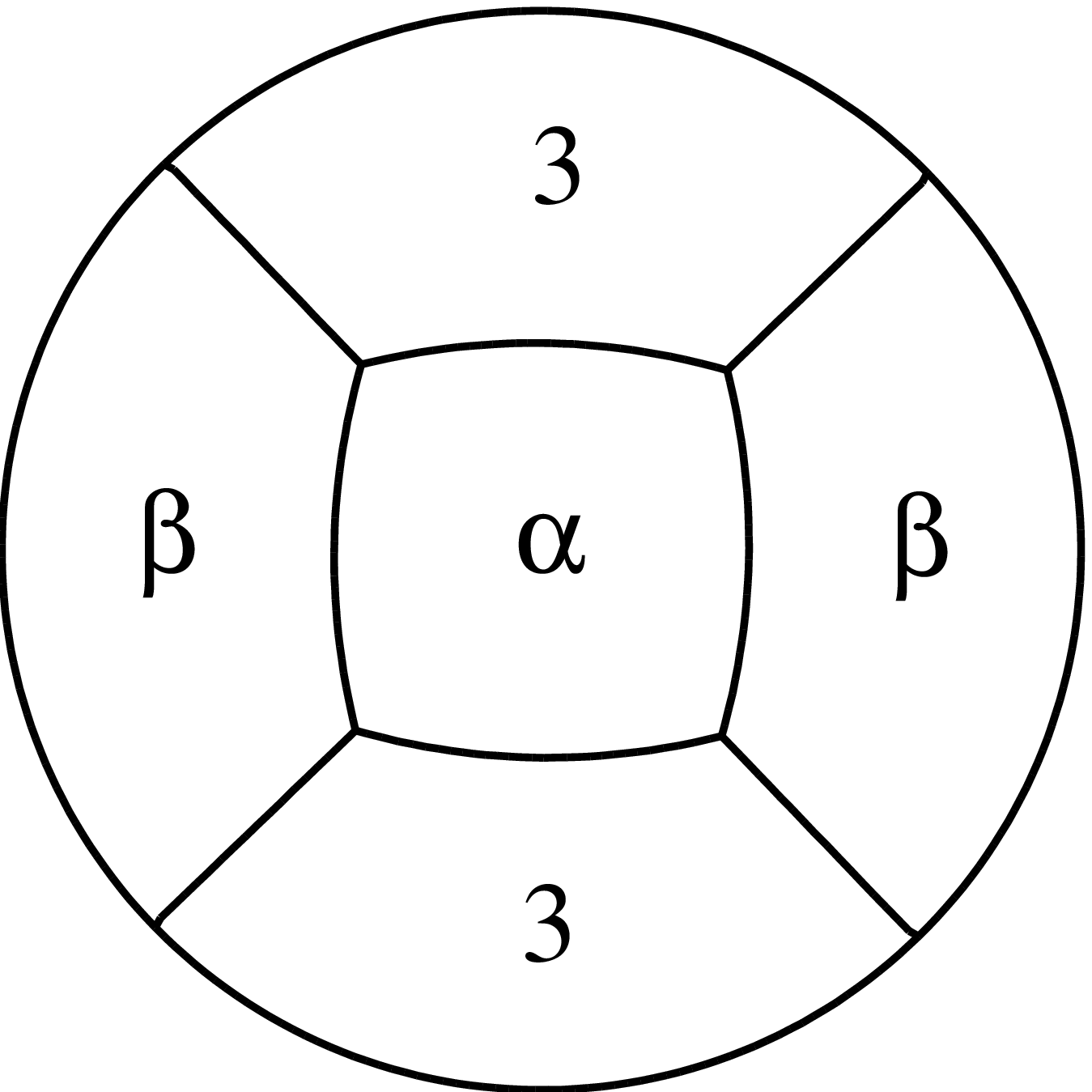}}
}
\\
\centering{
\subfigure[]{\label{conf6}\includegraphics[width=0.2\textwidth]{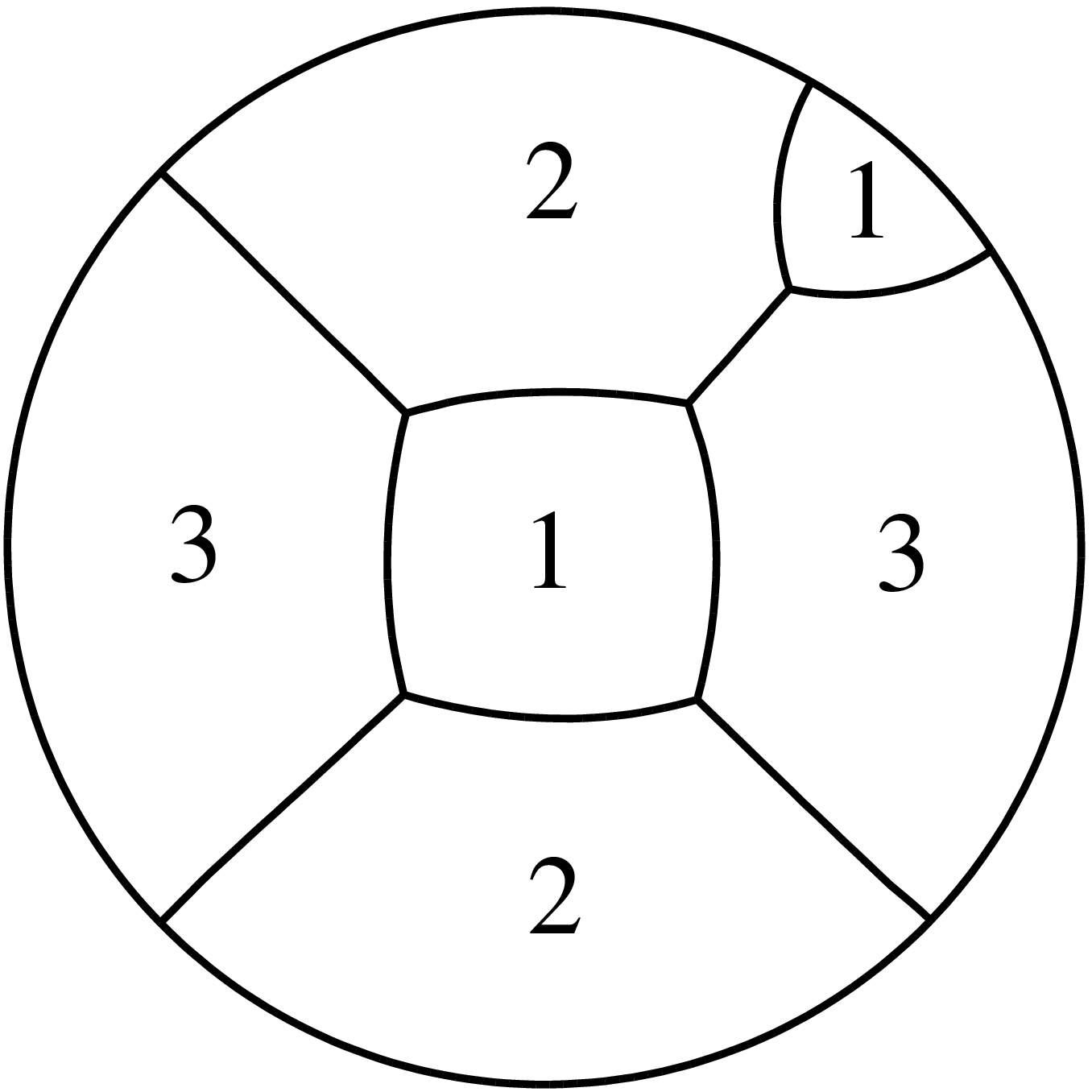}}
\hspace{0.1\textwidth}
\subfigure[$\alpha$, $\beta\in
\{2,3\}$]{\label{conf7}\includegraphics[width=0.2\textwidth]{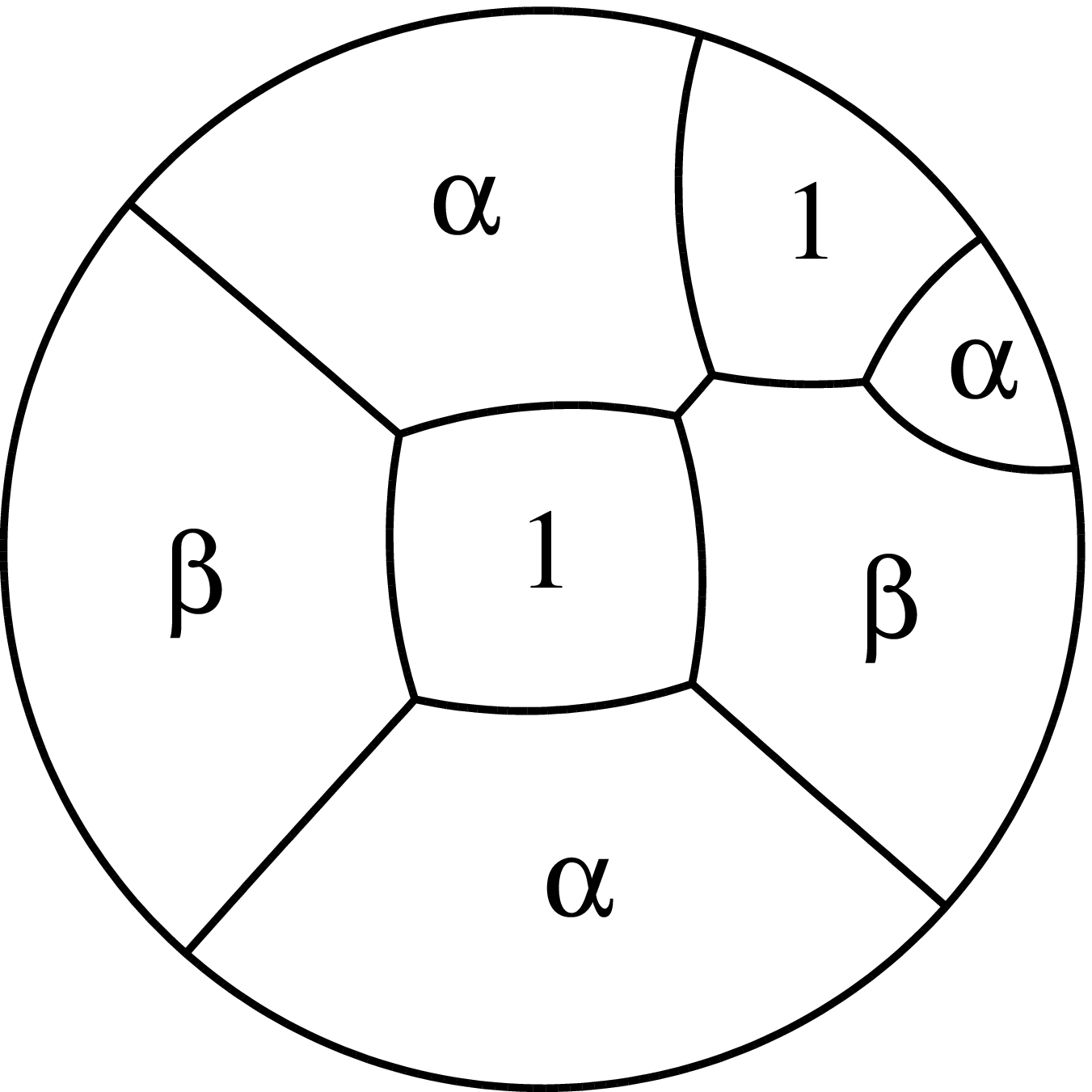}}
\hspace{0.1\textwidth}
\subfigure[]{\label{conf8}\includegraphics[width=0.2\textwidth]{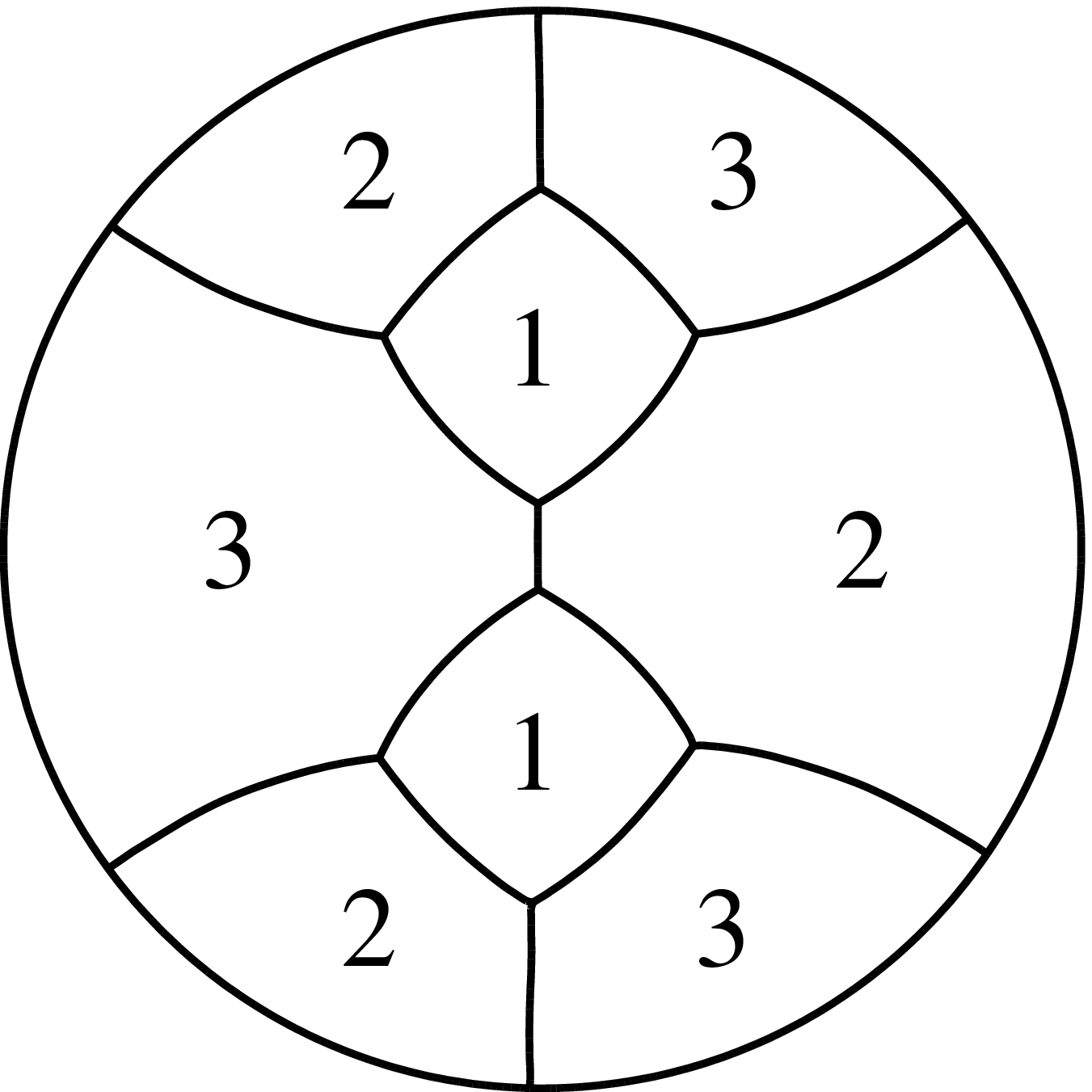}}
}
\\
\centering{
\subfigure[$\alpha$, $\beta\in
\{2,3\}$]{\label{conf9}\includegraphics[width=0.2\textwidth]{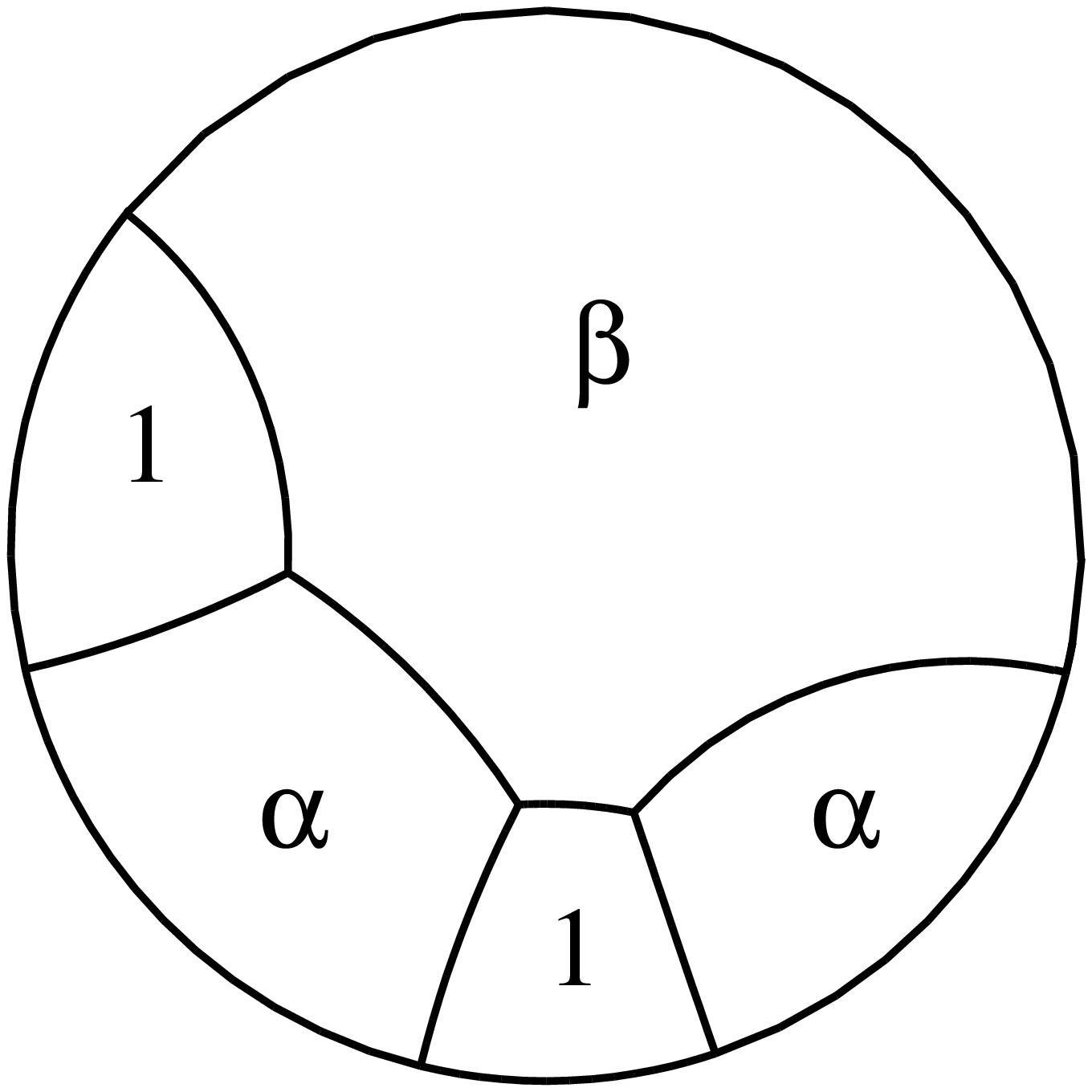}}
\hspace{0.1\textwidth}
\subfigure[]{\label{conf11}\includegraphics[width=0.2\textwidth]{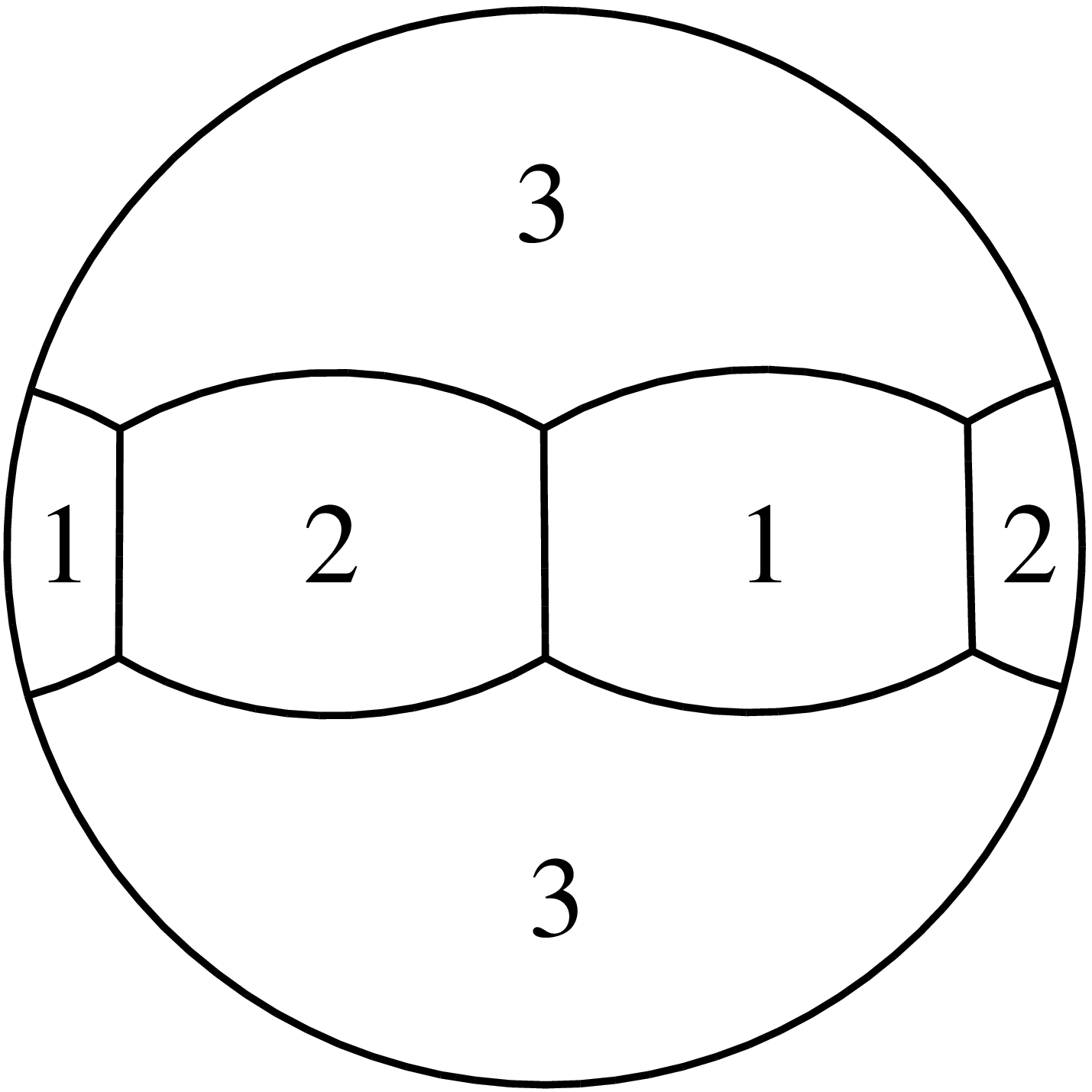}}
\hspace{0.1\textwidth}
\subfigure[]{\label{conf10}\includegraphics[width=0.2\textwidth]{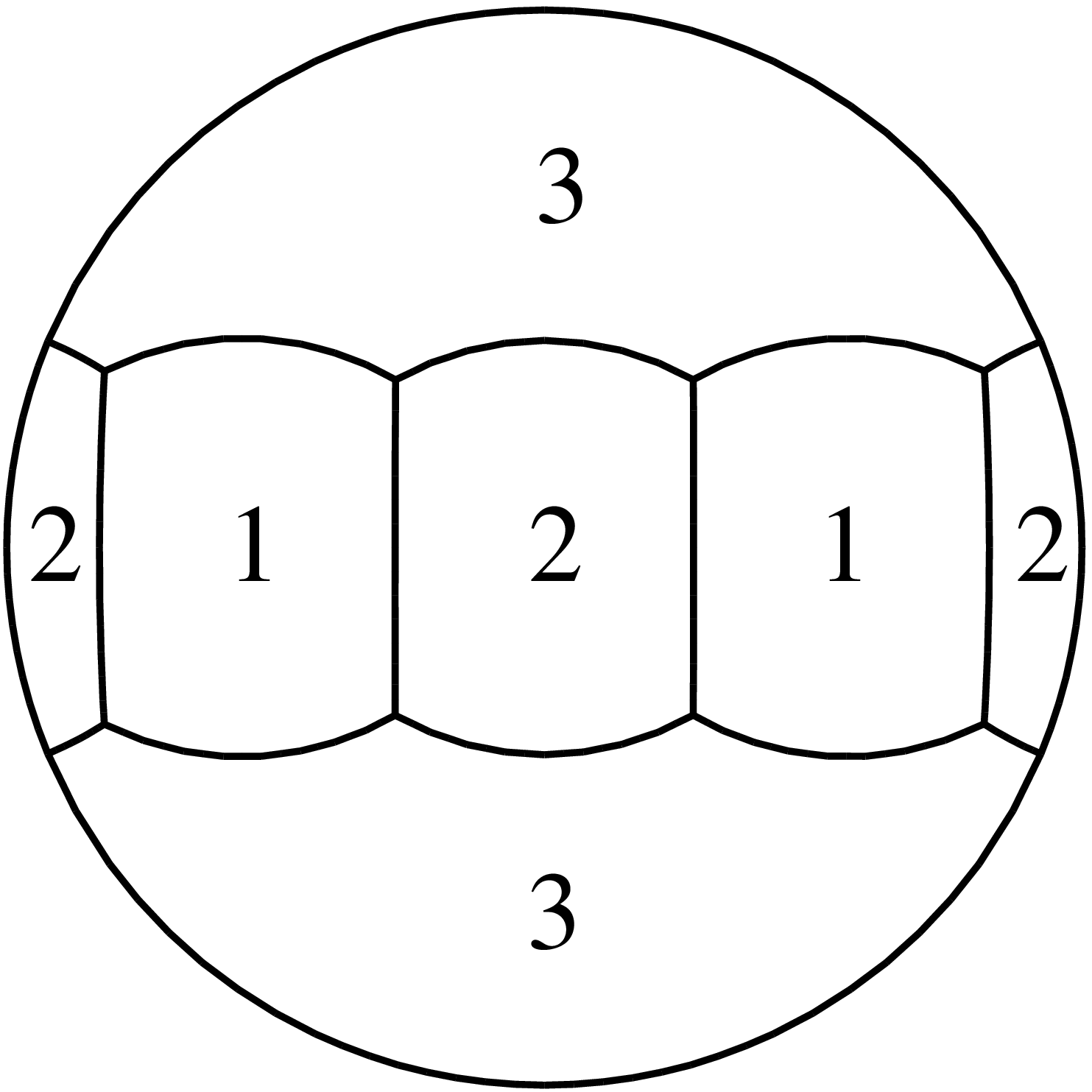}}
}
\\
\centering{
\subfigure[]{\label{conf12}\includegraphics[width=0.2\textwidth]{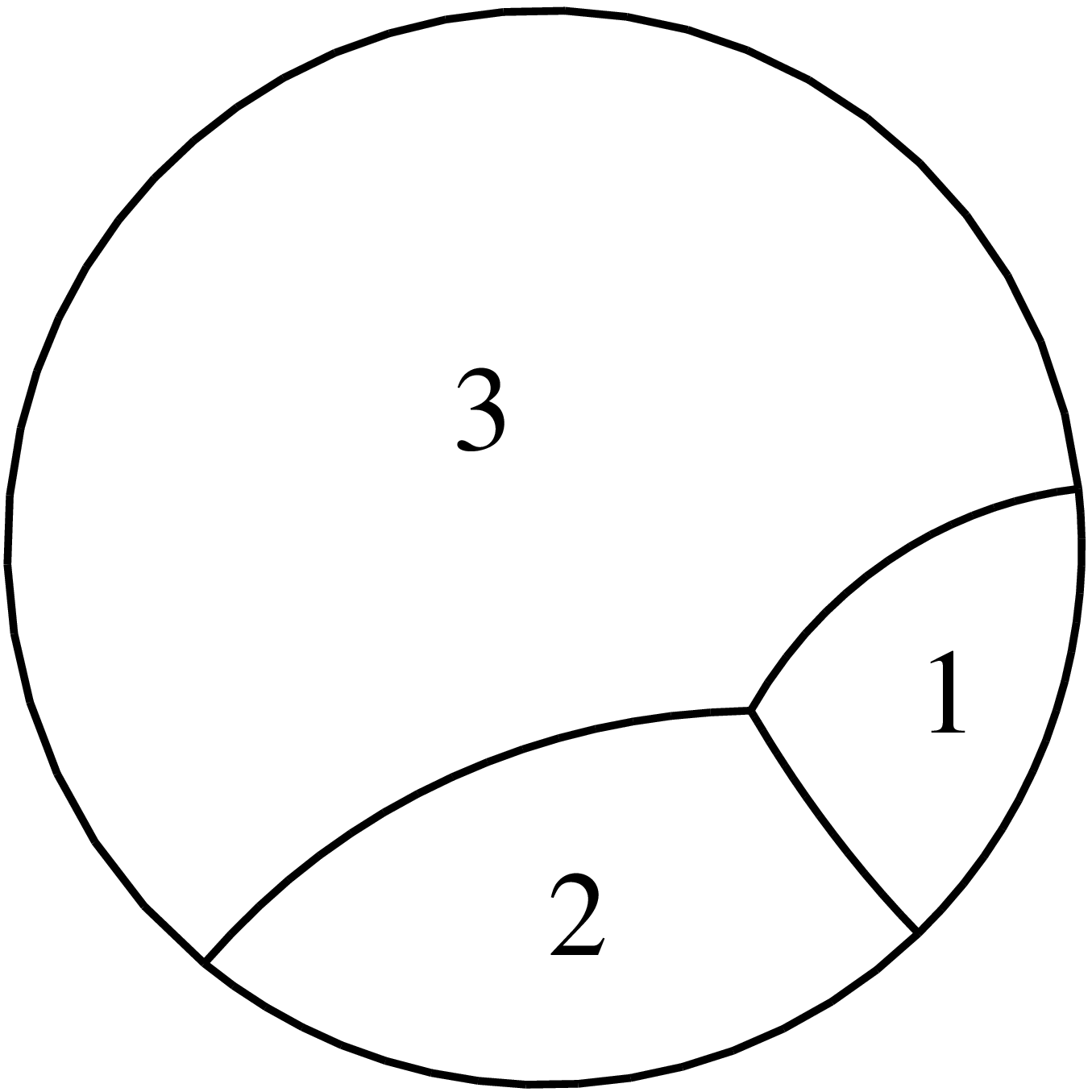}}
}
\caption{The ten possible configurations for minimizing graphs}
\label{fig:configs}
\end{figure}

Assume now that $p_1\ge p_2\ge p_3$, with $p_1>p_3$.  Then $R_1$
cannot have hexagonal components and so has at most two components, by
Lemma~\ref{le:cota}.  An interior component of $R_1$ has an even
number of edges in its boundary and cannot be a $2$-component by
Lemma~\ref{le:2-component}.  So any interior component of $R_1$ is a
quadrilateral.  A boundary component of $R_1$ will have by
Gauss-Bonnet Theorem three or four edges in its boundary.

Suppose $R_1$ is connected.  If it touches $\ptl D$ and has three
edges, we have the standard configuration \ref{conf12} and if it has
four edges, we have configuration \ref{conf1}.  If $R_1$ is interior ,
it will have only four edges, corresponding to configuration
\ref{conf4}.

Suppose now $R_1$ has two connected components, $A$ and $B$.  They can
be interior or boun\-da\-ry components.  We study each case.

If both are boundary components, as before, they can have three or
four edges.  If $A$ and $B$ have three edges, by the connectedness of
$C$, the only possibility is configuration \ref{conf1}.  If $A$ has
four edges and $B$ has three edges, the only possibility is
configuration \ref{conf9}; and if $A$ and $B$ have four edges, we will
have configurations \ref{conf3} and \ref{conf4}.

If $A$ is an exterior component and $B$ is an interior component, as
above, $A$ will have three or four edges, and $B$ will have four.  In
the first case, we will get configuration \ref{conf6}, and in the
second one, configurations \ref{conf7} and \ref{conf11}.

If both components are interior, they will have four edges.  A
component of the region with the smallest pressure cannot be interior
by Gauss-Bonnet Theorem.  Then the only possibilities are
configurations \ref{conf8} and \ref{conf10}.
\end{proof}

\section{Unstable and non-minimizing configurations}
\label{sec:unstable}

\begin{lemma}
\label{le:extension}
Let $C_{12}\subset D$ be a circle or segment meeting $\ptl D$
orthogonally between two regions $R_{1}'$, $R_{2}'$ with associated
pressures $p_{1}$ and $p_{2}$ $($the geodesic curvature of $C_{12}$
w.~r.~t.~the normal pointing into $R_{1}'$ equals $p_{1}-p_{2}$$)$.

Then, given $v\in C_{12}$, there exist unique curves $C_{23}$,
$C_{31}$ with constant geodesic curvature yielding a standard graph. 
Moreover, the pressure $p_{3}$ of the third region is a monotone
function of $v$ taking on all values between $-\infty$ and $+\infty$,
and the regions $R_{3}^v$ are nested.
\end{lemma}

\begin{proof}
Fix $v\in C_{12}$.  Let $q$ be a point in $C_{12}\cap\ptl D$. 
Consider the M\"obius transformation
\[
f(z)=\frac{i(z+q)}{q-z},
\]
which takes the disk $D$ to the upper half-plane and sends $q$ to
infinity.  Then $f(C_{12})$ is a straight line $L$.  Assume that there
are $C_{23}$, $C_{31}$ curves with constant geodesic curvature meeting
$C_{12}$ at $v$ so that $C_{12}\cup C_{23}\cup C_{31}$ is standard. 
Since $f$ is conformal, the sum of the geodesic curvatures of
$f(C_{ij})$ is zero.  Moreover $f(C_{23})$, $f(C_{31})$ intersect the
real axis orthogonally.  It is not difficult to see that $f(C_{23})$,
$f(C_{31})$ are circles centered at the real axis, with the same radius
by the balancing condition~\eqref{eq:cocycle}.  Thus $f(C_{23})$,
$f(C_{31})$ are unique and so are $C_{23}$ and $C_{31}$.  It is clear
that the regions determined by $f(C_{23})$ and $f(C_{31})$ are 
nested, which implies that $R_{3}^{v}$ are nested.

Let $d$ be the distance from $f(v)$ to the real axis.  By applying the
inverse of $f$, it is possible to compute the geodesic curvatures
$h_{31}$ and $h_{32}$ in terms of $d$, obtaining
\begin{align}
\label{eq:moeb1}
h_{31}&=\frac{1}{4}\,\bigg(-\sqrt{3}\,d - 2\,x + 
\frac{\sqrt{3}\,(1+x^2)}{d}\bigg),
\\
\label{eq:moeb2}
h_{32}&=\frac{1}{4}\,\bigg(-\sqrt{3}\,d + 2\,x + 
\frac{\sqrt{3}\,(1+x^2)}{d}\bigg).
\end{align}

Hence $p_3=p_1+h_{31}$, which decreases from $+\infty$ to $-\infty$
when $d$ moves from $0$ to $+\infty$.
\end{proof}

\begin{proposition}
\label{prop:extension}
Given three pressures, there is a standard graph separating the disk
into three regions with the given pressures.  Such a graph is unique
up to a rigid motion of the disk.
\end{proposition}

\begin{proof}
Given two pressures $p_{1}$, $p_{2}$, there is a circle or segment
$C_{12}$ separating $D$ into two regions with constant geodesic
curvature $p_{1}-p_{2}$ meeting $\ptl D$ orthogonally.  $C_{12}$ is
{\em unique} up to a rigid motion of the disk.  Using
Lemma~\ref{le:extension} we can find a vertex $v\in C_{12}$ and {\em
unique} curves $C_{23}$, $C_{31}$ so that $C_{12}\cup C_{23}\cup
C_{31}$ is a standard graph separating the disk into three regions
with the given pressures.  Uniqueness follows from the construction.
\end{proof}

\begin{proposition}
\label{prop:extension-2}
Let $C_{1}$, $C_{2}$, $C_{3}$ be circles or lines meeting at $120$
degree angles at some interior point of $D$, satisfying the balancing
condition~\eqref{eq:cocycle}.  If $C_{1}$ and $C_{2}$ meet $\ptl D$
orthogonally, then so it does $C_3$.
\end{proposition}

\begin{proof}
Let $\Om$ be the region enclosed by $C_{1}$, $C_{2}$ and $\ptl D$. 
Apply the disk onto the upper half-plane by means of a M\"obius map
$f$ sending $p\in\ptl D-\ptl\Om$ to infinity.  It is enough to show
that $f(C_{3})$ meets orthogonally the real axis.  As in the proof of
Lemma~\ref{le:extension}, $f(C_{1})$ and $f(C_{2})$ meet themselves at
$120$ degrees and the real axis at $90$ degrees.  As $f$ is a M\"obius
transformation, the balancing condition~\eqref{eq:cocycle} is
preserved.  In case $f(C_{3})$ is a line, it meets the real axis
orthogonally.  If $f(C_{3})$ is a circle, then $f(C_{1})\cup
f(C_{2})\cup f(C_{3})$ is a standard planar double bubble, for which
the centers of the circles are known to be aligned.  The centers of
$f(C_{1})$ and $f(C_{2})$ lie in the real axis, and hence also the
center of $f(C_{3})$.  So we conclude that $f(C_{3})$ meets the real
axis orthogonally.
\end{proof}

The proof of Lemma~\ref{le:extension} establishes the existence of a
deformation of one of the regions along one of the edges.  More
precisely we have

\begin{proposition}
\label{prop:moebius1}
Given a stationary graph $C$ with a boundary 3-component $\Om$, there
exists a variation of $C$ that
\begin{itemize}
\item[(i)] strictly increases the pressure of $\Om$ while keeping the
other pressures unchanged, and
\item[(ii)] strictly decreases the area of $\Om$, and
\item[(iii)] leaves invariant the edges of $C$ not lying in $\ptl\Om$.
\end{itemize}
\end{proposition}

\begin{proposition}
\label{prop:moebius2}
Let $C$ be a stationary graph in which a region has two boundary
3-compo\-nents.  Then $C$ is unstable.
\end{proposition}

\begin{proof}
Let $\Om_1$, $\Om_2$ be boundary 3-components of the same region.  On
each $\Om_i$, consider the variation given by Proposition
\ref{prop:moebius1}.  The normal components $u_i$ of the associated
variational vector fields have disjoint supports and satisfy
$Q(u_i,u_i)<0$.  By Propositions~\ref{prop:extension} and
\ref{prop:extension-2}, $\Om_1$ and $\Om_2$ are congruent so that
$u=u_1 - u_2$ satisfies the mean value conditions \eqref{eq:dareai}. 
Hence the graph is unstable.
\end{proof}

\begin{theorem}
\label{th:areas}
Given three areas $a_{1}$, $a_{2}$, $a_{3}$ such that
$a_{1}+a_{2}+a_{3}=\text{\em area}(D)$, there is a~unique standard
graph, up to rigid motions of the disk, separating $D$ into three
regions of areas $a_{i}$.
\end{theorem}

\begin{proof}
Consider two standard graphs: $C=C_{12}\cup C_{23}\cup C_{31}$, and
$C'=C'_{12}\cup C'_{23}\cup C'_{31}$ so that 
$a_{i}=\text{area}(R_{i})=\text{area}(R_{i}')$, $i=1$, $2$, $3$.

In case $h_{12}=h_{12}'$ we can apply a rigid motion of the disk to 
$C$ until $C_{12}$ and $C_{12}'$ coincide near $\ptl D$. As the areas 
of the enclosed regions are equal, Lemma~\ref{le:extension} implies 
that $C=C'$.

Assume that $h_{12}>h_{12}'$.  By Lemma~\ref{le:extension} we can
continuously decrease the pressure $p_{1}$ (while keeping constant
$p_{2}$ and $p_{3}$) until we get another standard graph $C''$ with
$h_{12}''=h_{12}$.  For this new graph $C''$ we get
$\text{area}(R_{1}'')>a_{1}$, $\text{area}(R_{i}'')<a_{i}$, for $i=2$,
$3$.  Now we can apply an isometry of the disk to $C''$ so that
$C_{12}''$ and $C_{12}'$ coincide near $\ptl D$.  As
$\text{area}(R_{1}'')>a_{1}$, Lemma~\ref{le:extension} implies that
$\text{area}(R_{2}'')>a_{2}$, which gives us a contradiction.
\end{proof}

Given two edges in a graph $C$, we will say that they are {\em
cocircular} if they have the same center.  A {\em cocircular}
4-component will be a 4-component with two cocircular opposite edges.

\begin{remark}
\label{re:sequence}
It is easy to check that in a sequence of interior 4-components, if
any of them is cocircular, then all the 4-components are.
\end{remark}

\begin{lemma}[{\cite[Lemma~5.38]{W}}]
\label{le:movement}
Suppose we have a stationary graph with a sequence of at least three
cocircular 4-components, so that the first and the last are boundary
components and the remaining are interior ones $($cocircularity refers
to the edges of the boundary components meeting $\ptl D)$.  Assume
further that the components out of the chain belong to the same
region.  Then there is a continuous movement preserving perimeter and
areas which creates an illegal meeting, so that the graph is not
minimizing.
\end{lemma}

\begin{proof}
Order the 4-components so that $\Om_1$ and $\Om_n$ are the boundary
ones, and $\Om_i$ meets $\Om_{i+1}$.  Let $c_i$ be the center of the
cocircular arcs of $\Om_i$.  We can move these points, without
changing neither $d(c_i,c_{i+1})$ nor $d(c_1,0)$ and $d(c_n,0)$ in
such a way that $c_1$ and $c_n$ get closer and closer.  With this
movement of the points, we obtain a deformation of the graph which
preserves perimeter and the areas of the regions and will create an
irregular meeting.  Hence the graph cannot be minimizing.
\end{proof}

\begin{lemma}[{\cite[Lemma~5.3]{W}}]
\label{le:simetria}
Let e, f and g be three consecutive edges of a component, and let
$v_1$ and $v_2$ be the corresponding vertices.  Suppose e and g have
the same geodesic curvature, and the angles in each vertex are the
same.  Let R be the line of points equidistant from $v_1$ and $v_2$.

Then e and g are interchanged by the symmetry about R. Moreover, if e
and g are cocircular, the common center lies in R, and if e and g are
not cocircular, R coincides with the line of points equidistant from
the centers of e and g.
\end{lemma}

Given a graph $C$ and a Jacobi function $u$ defined on it, we will say
that a point $x$ in $C$ is a {\em nodal point} if $u(x)=0$.  A {\em
nodal region} will be a connected component of the complementary in
$C$ of the set of nodal points.

\begin{proposition}[{\cite[Proposition~5.2]{HMRR}}]
\label{prop:nodalregions}
Let $C$ be a stationary graph separating the disk into three regions. 
Assume that there exists a Jacobi function with at least four nodal
regions such that the nodal points are not vertices of the graph. 
Then C is unstable.
\end{proposition}
 
\begin{proof}
Let $u$ be the Jacobi function, and $N_1, \ldots, N_4$ nodal regions. 
Assume that the graph $C$ is stable.  For $i=1,\ldots,4$ define $u_i$
as the restrictions of $u$ to $N_i$ extended by zero to the whole
graph.  It is possible to find a nontrivial linear combination $v$ of
$u_1$, $u_2$, $u_3$ so that the mean value conditions
\eqref{eq:dareai} are satisfied for $R_1$, $R_2$ and $R_3$.  Moreover,
$v$ is an admissible function and, by stability, is a Jacobi function. 
As $v$ vanishes on a subset of $C$ containing $N_4$ and has nontrivial
support contained in $N_1\cup N_2\cup N_3$, there is an edge $\ell$ so
that $v$ vanishes on an open interval of $\ell$ but it is not
identically zero on $\ell$.  As $v$ is a Jacobi function (a solution
of a second order o.d.e.) this gives us a contradiction.
\end{proof}

Now we discard configuration~\ref{conf11} by a geometrical argument. 
The reader may compare this result with \cite[Lemma~5.2.10]{corneli}

\begin{proposition}
\label{prop:ultimaconf}
Configuration $\ref{conf11}$ is not minimizing.
\end{proposition}

\begin{figure}[h]
\centering{\includegraphics[width=0.7\textwidth]{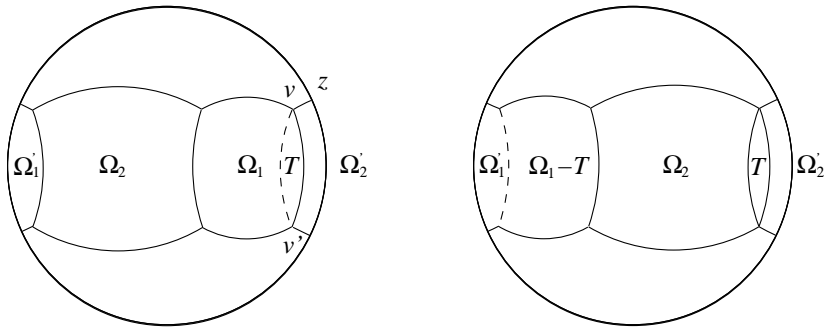}}
\caption{Deformation of a chain
 of symmetric 4-components obtaining a graph with irregular vertices}
\label{fig:conf11}
\end{figure}

\begin{proof}
Suppose this configuration is minimizing.  We will denote by
$\Om_{1}\subset R_1$, $\Om_{2}\subset R_2$ the interior components,
and by $\Om_{1}'\subset R_1$, $\Om_{2}'\subset R_2$ the boundary
components.

Consider any component of $R_1$ or $R_2$ and suppose the edges
separating such a component and $R_3$ are cocircular.  By
Remark~\ref{re:sequence} we will have a sequence of four cocircular
4-components and using Lemma~\ref{le:movement} the configuration is
not minimizing.

Hence the considered edges cannot be cocircular, and applying
Lemma~\ref{le:simetria}, it is easy to check that there exists a
horizontal symmetry of the chain of 4-components which is in fact a
symmetry of the whole configuration.

In each interior component, using Lemma~\ref{le:simetria} and taking
into account that any pair of opposite edges cannot be cocircular, it
can be seen that there exists a vertical symmetry, orthogonal to the
horizontal one.

If $p_1 = p_2$, we have that $\Om_1$ and $\Om_2$ are congruent, so
that they have the same area.  Hence, we can interchange them, and
after eliminating unnecessary edges, we obtain a new configuration
enclosing the same areas with strictly less perimeter, which gives a
contradiction.

If $p_2 = p_3$, applying Gauss-Bonnet Theorem to $\Om_{2}$ we obtain a
contradiction.

So we can suppose now $p_1>p_2>p_3$. Let us distinguish two cases:

Let $z$ be the upper point in $\ptl\Om_{2}'\cap\ptl D$.  Let $\theta$
be the angle between the line $0z$ and a horizontal line.  Assume that
$\theta$ is greater than or equal to $\pi/4$.  Then, the second
coordinate of $z$ will be greater than or equal to
$\cos{\pi/4}=\sqrt{2}/2$.  We can see $\ptl \Om_{2}'$ as a vertical
graph, so the length of $\ptl \Om_{2}'$ will be greater than or equal
to $\sqrt{2}$.  By similar arguments, taking into account that the
point of $\ptl \Om_{2}$ with maximum second coordinate will be higher
than $z$, the length of $\ptl \Om_{2}$ will be greater than or equal
to $2\sqrt{2}$.  Both quantities add up to more than $3$, so by
Lemma~\ref{le:cotaperfil}, this configuration cannot be minimizing.

Assume now that $\theta$ is less than $\pi/4$.  Let $l$ be the edge
separating $\Om_{1}$ and $\Om_{2}'$, and $v$, $v'$ its vertices. 
Consider a new edge $\widetilde{l}$, the symmetric of $l$ about the
segment $\overline{vv'}$, and let $T$ be the region enclosed between
$l$ and $\widetilde{l}$.  Move $\Om_2$ in the horizontal direction
until one of its edges coincides with $\widetilde{l}$.  Reflect
$\Om_{1}-T$ about a vertical axis so that the reflection of
$\widetilde{l}$ coincides with the edge $\Om_2\cap \Om_{1}'$ (recall
that each interior 4-component has a vertical symmetry, so the side
edges have the same length and curvature).  After eliminating an
unnecessary edge, this new configuration will preserve length and
areas, but it will be irregular.  Hence, configuration \ref{conf11}
cannot be minimizing.  It only remains to check that $\Om_{1}$ and
$\Om_{2}$ will remain inside $D$ under this geometrical
transformation.

In order to prove this it is enough to show that the portion of the
original graph over the horizontal line $L$ passing through $z$ stays
inside $D$.  The portion $\ell_2$ of the upper edge of $\Om_2$ over
$L$ makes an angle $\theta$ with $L$, and can be translated
horizontally to touch $z$ at its boundary.  Since $\theta<\pi/4$,
$\ell_2$ must lie inside $D$.  This implies that the transformation of
$\Om_2$ stays inside $D$.  The upper edge $\ell_1$ of $\Om_1$ has
larger geodesic curvature than the upper edge of $\Om_2$ and makes a
smaller angle with the horizontal line passing through $v$.  This is
enough to conclude that $\ell_1$ must lie inside $D$.
\end{proof}

\section{Proof of the theorem} 
\label{sec:teorema}

In this section we prove the main theorem.

\begin{theorem}
\label{te:main} 
Let $C\subset D$ be a minimizing graph for three given areas.  Then
$C$ is a unique standard graph.
\end{theorem}

\begin{proof}
By Lemma~\ref{le:possibilities} the graph $C$ must be one of the
listed in Figure~\ref{fig:configs}.

Configurations \ref{conf1} and \ref{conf3} are unstable by
Proposition~\ref{prop:moebius2}.

Configuration \ref{conf4} is also unstable: it is easy to check that
the edges of $C_{12}$ are not cocircular, so that by
Remark~\ref{re:sequence} and Lemma~\ref{le:simetria}, we have a
vertical symmetry in this configuration.  By Lemma~\ref{le:movement},
the edges of $C_{13}$ cannot be cocircular if $C$ is minimizing. 
Applying again Lemma~\ref{le:simetria} we get an horizontal symmetry. 
Both axes of symmetry will meet orthogonally at the origin, so we can
consider the Killing field generated by rotations about $0$, which
yields a Jacobi function $u$ vanishing on four points, one on the
interior of each edge of the central 4-component.  So $u$ has at least
four nodal regions and we conclude by
Proposition~\ref{prop:nodalregions} that $C$ is unstable.

We now eliminate configurations \ref{conf6} and \ref{conf7}.  These
configurations present an interior 4-component of $R_1$, with three
incident edges meeting the exterior of the disk.  If we extend the
fourth edge, it will meet $\ptl D$ orthogonally by
Proposition~\ref{prop:extension-2}, and we will obtain a configuration
of type \ref{conf4}.  Hence, the interior 4-component has two
orthogonal symmetries meeting at the origin, and we conclude as before
the existence of four nodal regions.

Consider now configuration \ref{conf8}.  Fix an interior 4-component
$\Om$ of $R_1$.  If we extend the edge leaving the boundary of $\Om$
that does not reach $\ptl D$, it will meet $\ptl D$ orthogonally due
to the existence of a symmetry of $\Om$ which is in fact a symmetry of
the disk.  In this way we obtain a configuration of type \ref{conf4}. 
As above, $\Om$ will have two orthogonal symmetries meeting at $0$ and
so we can get four nodal regions yielding instability.

Consider now configuration \ref{conf9}.  Applying
Lemma~\ref{le:simetria}, the two 4-components will be symmetric about
two lines $r_1,\,r_2$ passing through the center of the disk (the
corresponding edges are not cocircular).  Let $q_1,\ q_2$ be the
intersection points of each line with the interior edges of these
components, that will be zeros of the Jacobi field $u$ induced by the
one-parameter group of rotations about the origin.  The reflection of
$q_1$ with respect to $r_2$ lies in the boundary of the 3-component of
$R_1$ and it is not a vertex of the configuration.  This point is
clearly also a zero of $u$.  Then $u$ has four nodal regions and the
configuration is unstable.

By Proposition~\ref{prop:ultimaconf}, configuration \ref{conf11} is
not minimizing.

Configuration \ref{conf10} is unstable: if the top and bottom edges of
each component are cocircular then the configuration is not minimizing
by Lemma~\ref{le:movement}.  Otherwise we can find an horizontal
symmetry of the graph, which is also a symmetry of the disk by
Lemma~\ref{le:simetria}.  Each interior 4-component has a vertical
symmetry so that the interior components of $R_1$ are identical. 
Using the function equal to $+1$ on one of these components, equal to
$-1$ on the other component, and zero otherwise, we have obtained a
function satisfying the mean value conditions \eqref{eq:dareai} which
is negative for the index form.  So this configuration is unstable. 
We could also use the method of Proposition~\ref{prop:ultimaconf} to
see that this configuration is nonminimizing.
 
So the only remaining possibility is configuration \ref{conf12}, the
standard one.  Uniqueness for given areas comes from 
Theorem~\ref{th:areas}.
\end{proof}

\section{Final remarks}
\label{sec:final}
In this paper we have obtained that the problem of dividing the disk
into three areas has a unique solution in which all regions are
connected, as in the problem of partitioning the disk into two areas. 
It is natural to conjecture that

\begin{conjecture}
A minimizing graph separates the disk into connected regions. 
\end{conjecture}

If we consider the problem for $n$ regions, with $n\geq 4$, by
Lemma~\ref{le:cota} the region of largest pressure will have at most
$n-1$ nonhexagonal connected components and we can obtain by
combinatorial arguments a list of all possible minimizing
configurations.  Of course the number of candidates increases very
rapidly when the number of regions increases.  We believe that the
following conjectures are true

\begin{figure}[htp]
\centering{
\subfigure[]{\label{cuatro}\includegraphics[width=0.2\textwidth]{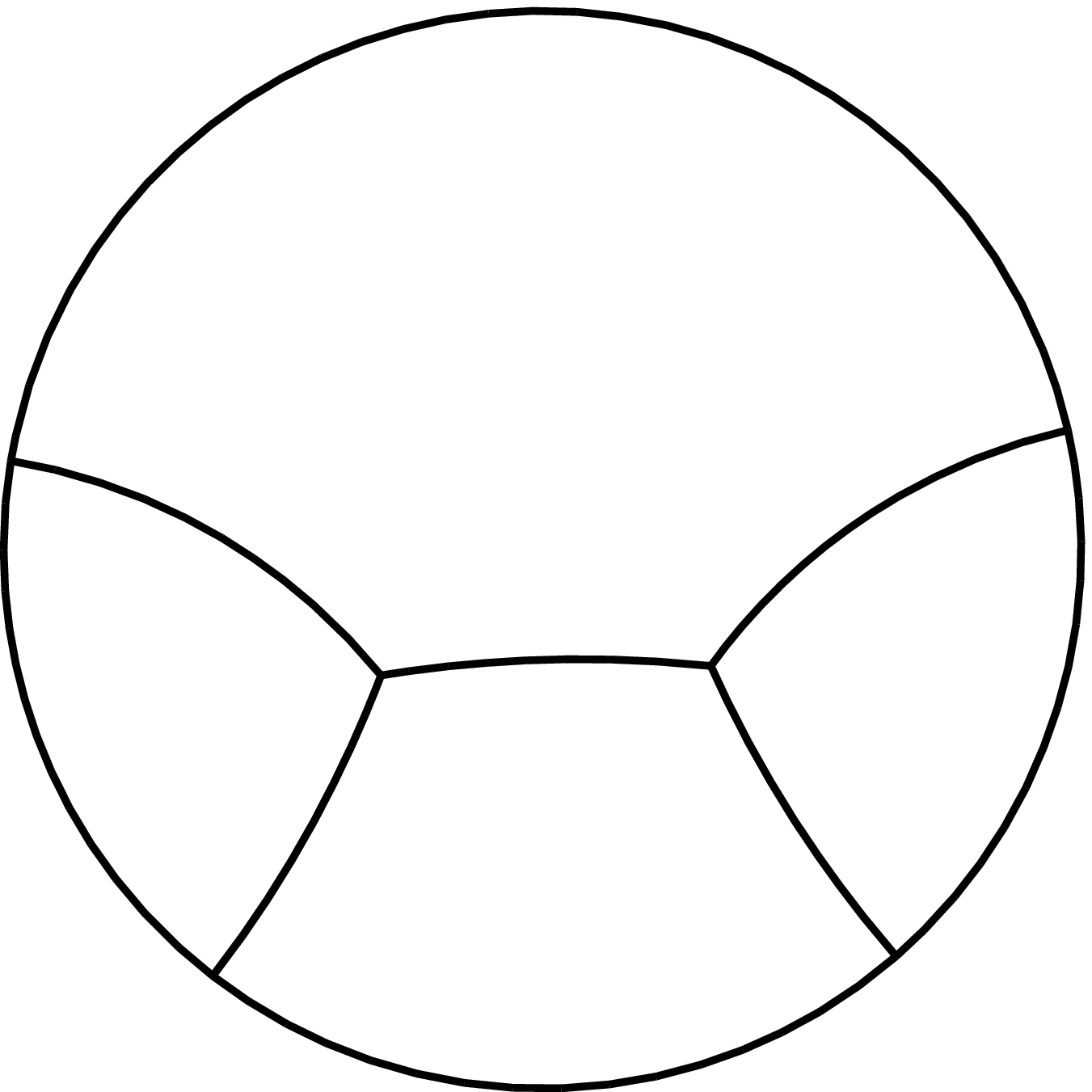}}
\hspace{0.1\textwidth}
\subfigure[]{\label{cinco}\includegraphics[width=0.2\textwidth]{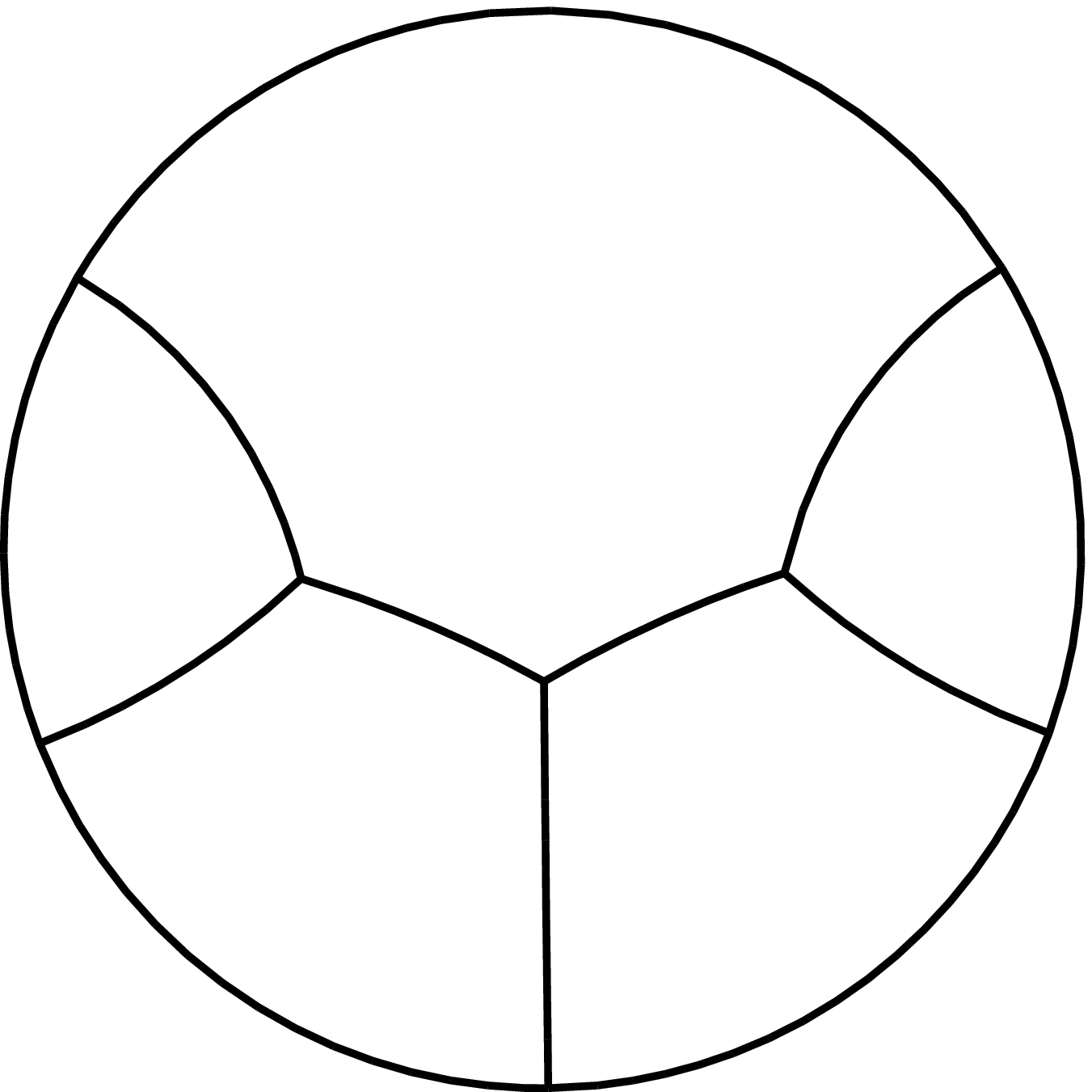}}
}
\caption{The conjectural configurations for $n=4$ and $n=5$}
\label{fig:conjectural}
\end{figure}

\begin{conjecture}
\label{conj:cuatro}
The least perimeter way of dividing the unit disk into four regions of
pres\-cribed areas is given by configuration \ref{cuatro} of
Figure~\ref{fig:conjectural}.
\end{conjecture}

\begin{conjecture}
\label{conj:cinco}
The least perimeter way of dividing the unit disk into five regions of
pres\-cribed areas is given by configuration \ref{cinco} of
Figure~\ref{fig:conjectural}.
\end{conjecture}

For each case, we believe that there is another possibly stable
configuration: for $n=4$, the configuration with three boundary
regions surrounding an interior one of three edges, and for $n=5$, the
one consisting in four boundary regions surrounding an interior region
of four edges.  But estimates we have done using Surface Evolver (Ken
Brakke, 1992) for equal areas show that they are nonminimizing. 
Furthermore, for $n=4$, if any of the areas tends to zero, we should
obtain in the limit the standard configuration for three areas, which
also discards the configuration described above at least for some
areas.  In the case $n=5$ we should have the same behaviour.

For $n=6$ we give the following conjecture

\begin{figure}[h]
\centering{\includegraphics[width=0.2\textwidth]{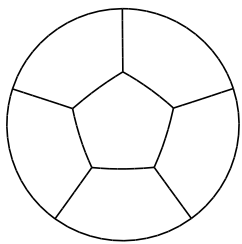}}
\caption{The conjectural configuration for $n=6$}
\label{fig:seisconj}
\end{figure}

\begin{conjecture}
\label{conj:seis}
The least perimeter way of dividing the unit disk into six regions of
prescribed areas is given by configuration of
Figure~\ref{fig:seisconj}.
\end{conjecture}

\begin{figure}[h]
\centering{
\subfigure[]{\label{seisuno}\includegraphics[width=0.2\textwidth]{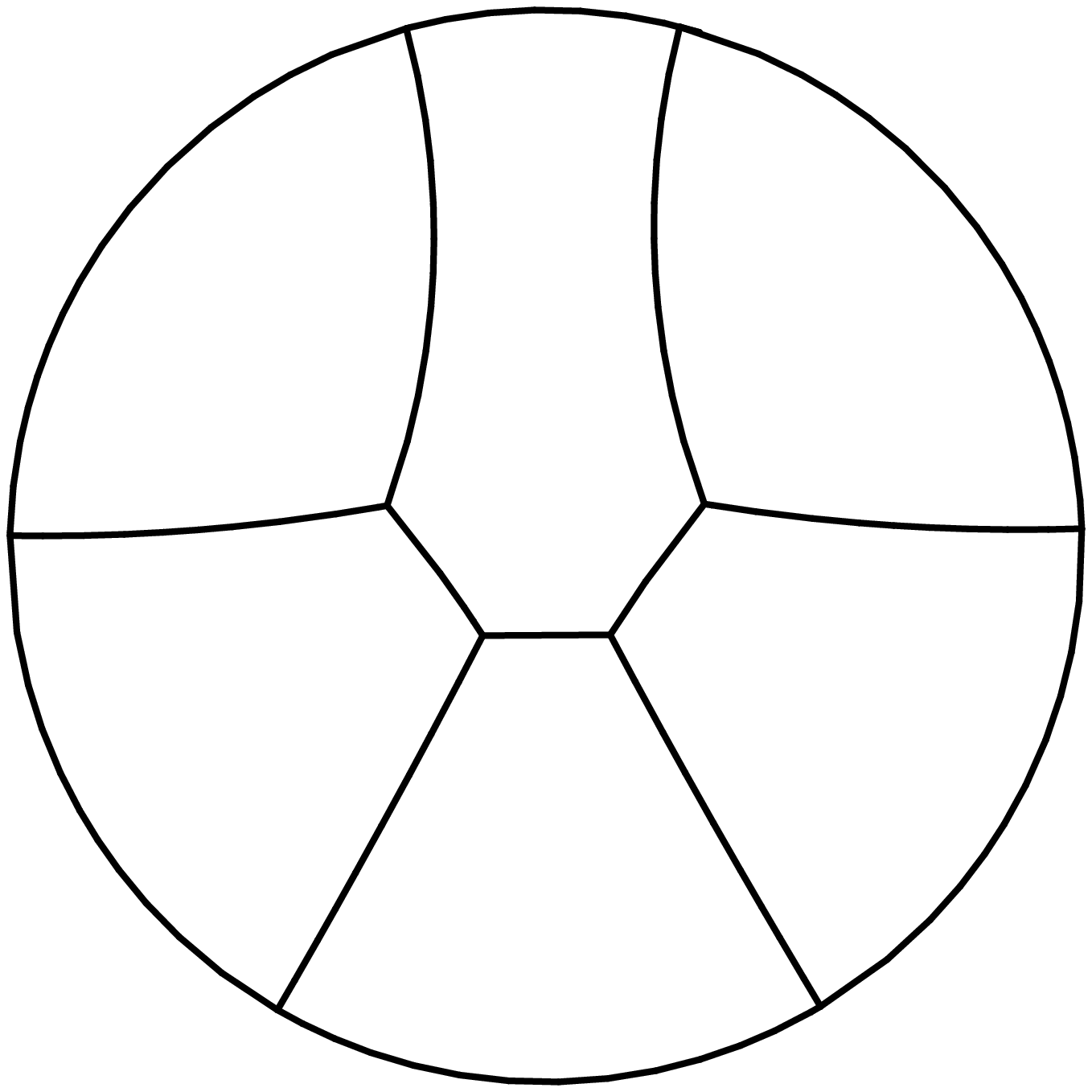}}
\hspace{0.05\textwidth}
\subfigure[]{\label{seisdos}\includegraphics[width=0.2\textwidth]{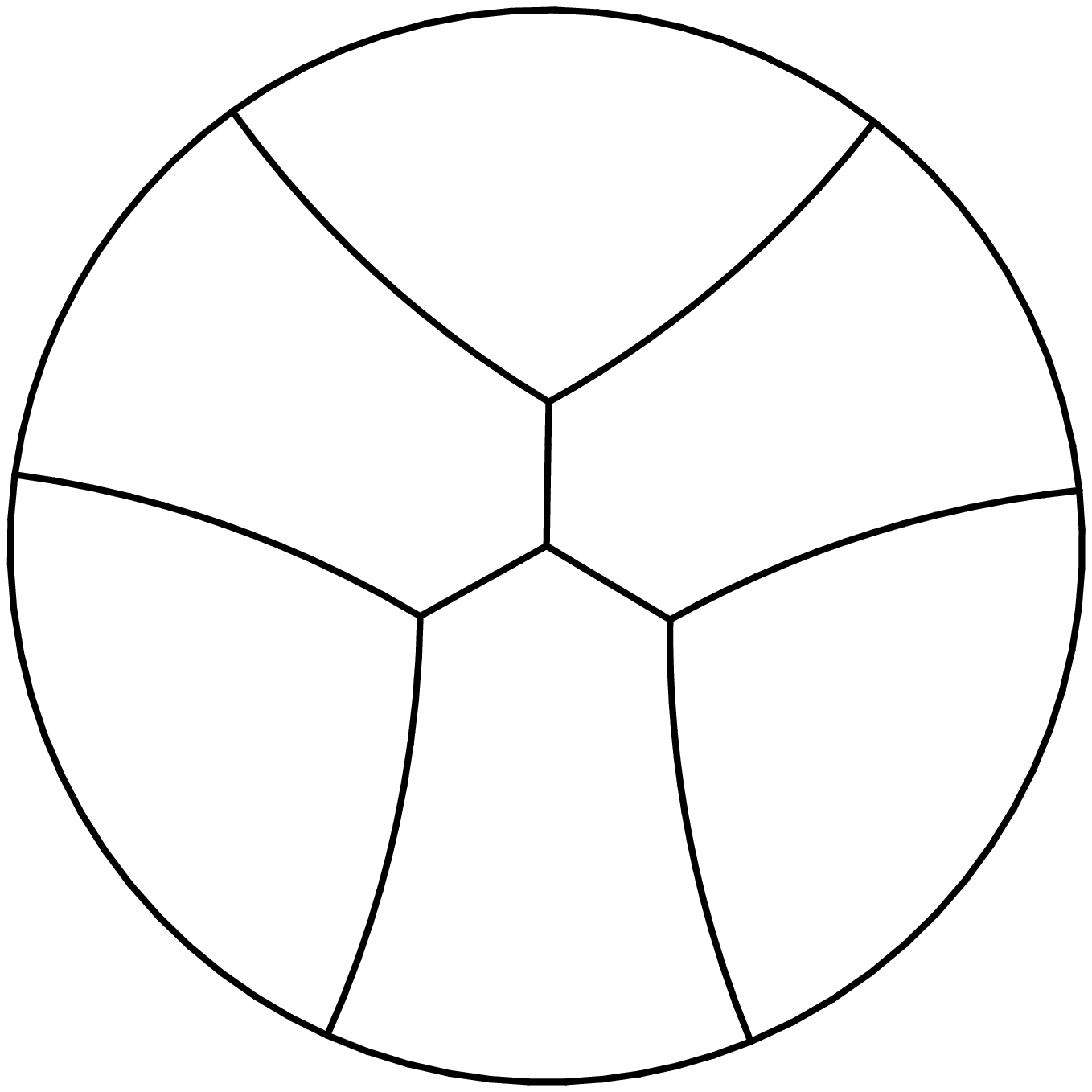}}
\hspace{0.05\textwidth}
\subfigure[]{\label{seistres}\includegraphics[width=0.2\textwidth]{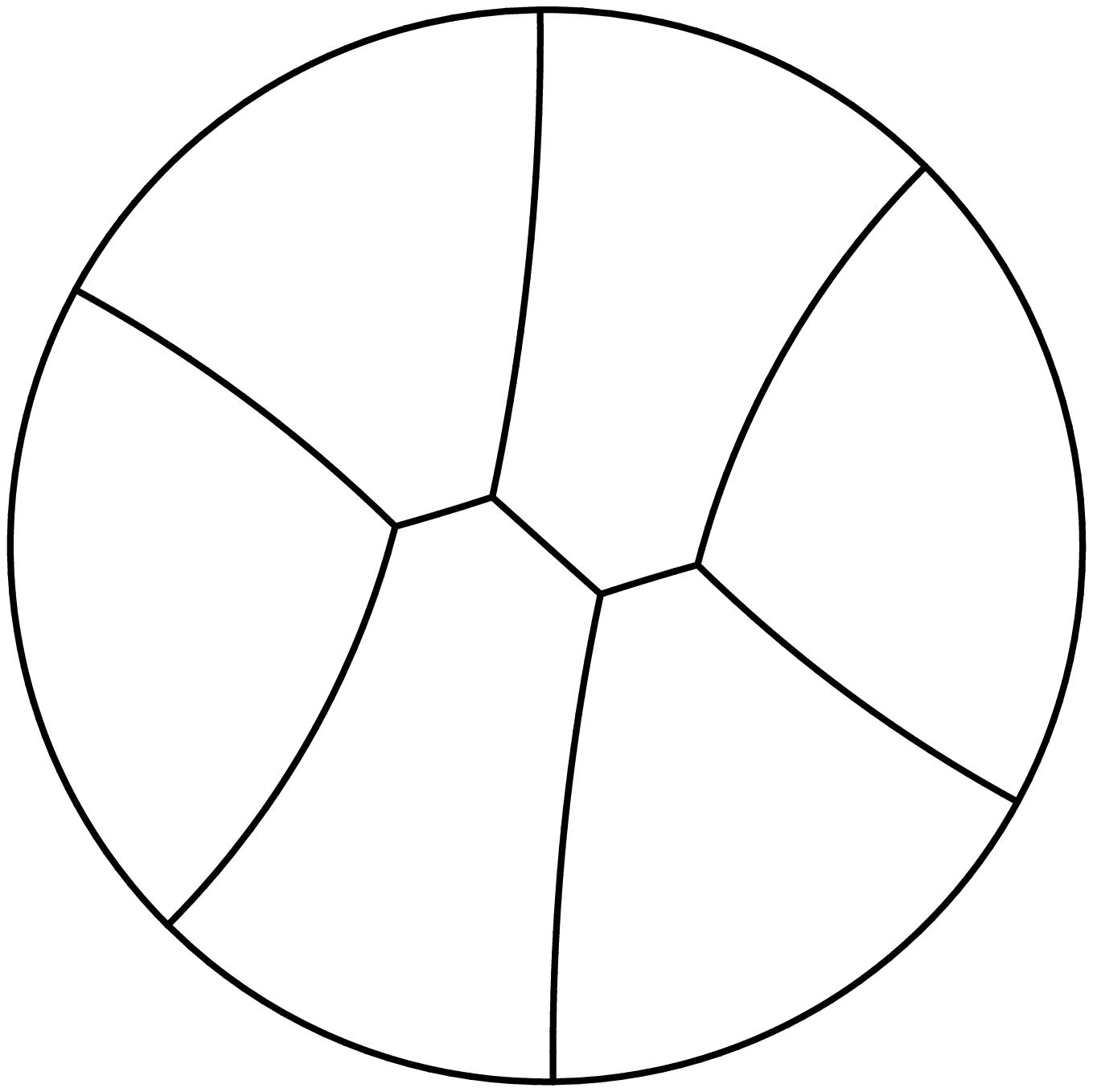}}
}
\caption{Some other configurations for $n=6$}
\label{fig:seisestables}
\end{figure}

As before, we belive that the configurations of
Figure~\ref{fig:seisestables} below are stable, but estimates done
with the Surface Evolver considering equal areas show that they are
nonminimizing.

\newpage

\providecommand{\bysame}{\leavevmode\hbox to3em{\hrulefill}\thinspace}
\providecommand{\MR}{\relax\ifhmode\unskip\space\fi MR }
\providecommand{\MRhref}[2]{%
  \href{http://www.ams.org/mathscinet-getitem?mr=#1}{#2}
}
\providecommand{\href}[2]{#2}

\end{document}